\documentclass[sort&compress,3p,fleqn]{elsarticle}
\usepackage{pdfpages,color,pdfsync}
\usepackage{graphicx}
\usepackage{amsmath,amssymb,amsthm}
\usepackage{hyperref}

\allowdisplaybreaks[4]

\numberwithin{equation}{section}

\newtheorem{Def}{Definition}[section]
\newtheorem{Thm}[Def]{Theorem}
\newtheorem{Lem}[Def]{Lemma}

\newtheorem{Coro}[Def]{Corollary}

\newcommand{\R}{\mathbb{R}}
\newcommand{\J}{\mathbb{J}}
\newcommand{\mC}{\mathbb{C}}
\newcommand{\Z}{\mathbb{Z}}

\newcommand{\T}{\mathbb{T}^2}

\newcommand{\W}{W_{\T}}

\newcommand{\I}{i}

\newcommand{\ud}{d}

\newcommand{\im}{\mathrm{Im}\,}

\newcommand{\supp}{\mathrm{supp}\,}
\newcommand{\Hess}{\mathrm{Hess}}

\newcommand{\bvec}[1]{\boldsymbol{#1}}

\newcommand{\vj}{\bvec{j}}

\newcommand{\vq}{\bvec{q}}
\newcommand{\vn}{\bvec{n}}
\newcommand{\vx}{\bvec{x}}
\newcommand{\va}{\bvec{a}}

\newcommand{\val}{\bvec{a}^0}
\newcommand{\vb}{\bvec{b}}

\newcommand{\p}{\partial}

\newcommand{\Div}{\nabla\cdot}

\newcommand{\la}{\left\langle}
\newcommand{\ra}{\right\rangle}
\newcommand{\nn}{\nonumber\\}

\newcommand{\SE}{Schr\"{o}dinger equation}

\newcommand{\wto}{\rightharpoonup}

\newcommand{\be}{\begin{equation}}
\newcommand{\ee}{\end{equation}}

\newcommand{\M}{W^{-1,1}}

\newcommand{\Bluemark}[1]{{#1}}
\newcommand{\blackmark}[1]{{#1}}

\begin{document}

\begin{frontmatter}
\title{Quantized vortex dynamics of the nonlinear Schr\"{o}dinger equation \blackmark{on torus} with non-vanishing momentum}

\author[Zhu]{Yongxing Zhu\corref{cor}}
\ead{zhu-yx18@mails.tsinghua.edu.cn}

\address[Zhu]{Department of Mathematical Sciences, Tsinghua University, Beijing, 100084, China}

\author[Bao]{Weizhu Bao}
\address[Bao]{Department of Mathematics, National University of Singapore, Singapore 119076, Singapore}
\ead{matbaowz@nus.edu.sg}

\author[Zhu]{Huaiyu Jian}
\ead{hjian@tsinghua.edu.cn}

\cortext[cor]{Corresponding author.}

\begin{abstract}
We derive rigorously the \blackmark{reduced dynamical law} for quantized vortex dynamics
of the nonlinear \SE{} on \blackmark{the torus} with non-vanishing momentum when the vortex core size
$\varepsilon\to0$. \blackmark{The reduced dynamical law is} governed by a Hamiltonian flow driven by a renormalized energy. A key ingredient is to construct a new canonical harmonic map to include the effect from the non-vanishing momentum into the dynamics. Finally, some properties of the reduced dynamical law
are discussed.
\end{abstract}


\date{}

\end{frontmatter}
\section{Introduction}\label{sec:intro}
In this paper, we study quantized vortex dynamics of the nonlinear \SE{} (NLSE) on \blackmark{the torus} \cite{RN8,Rubinstein1991}:
\begin{equation}\label{eq:maineq}
    \I \p_tu^\varepsilon(\vx,t)- \Delta u^\varepsilon(\vx,t)+\frac{1}{\varepsilon^2}\left(|u^\varepsilon(\vx,t)|^2-1 \right)u^\varepsilon(\vx,t)=0,\quad  \vx\in\T, \ t>0,
\end{equation}
with initial data
\begin{equation}\label{init}
u^\varepsilon(\vx,0)= u_0^\varepsilon(\vx),\quad \vx\in\T,
\end{equation}
where \blackmark{ $t$ is the time variable}, $\T=(\R/\Z)\times(\R/\Z)$ is the unit  torus, $\vx=(x,y)^T$ is the spatial coordinate,  $u^\varepsilon:=u^\varepsilon(\vx,t)$ is a complex-valued wave function or order parameter, $0<\varepsilon<<1$ is a dimensionless parameter which is used to characterize the core size of quantized vortices, and
$u_0^\varepsilon:=u_0^\varepsilon(\vx)$ is a given initial data.  It is well-known that the NLSE
\eqref{eq:maineq} conserves the {\sl mass} defined as \cite{LinXin1999}
\begin{equation}\label{eq:conserve M}
M(u^\varepsilon(t)):=\int_{\T} |u^\varepsilon(\vx,t)|^2d\vx\equiv \int_{\T} |u^\varepsilon_0(\vx)|^2d\vx:= M(u^\varepsilon_0), \qquad t\ge0;
\end{equation}
the {\sl momentum} defined as \cite{RN8,RN23,LinXin1999}
\begin{equation}\label{eq:conserve Q}
\bvec{Q}(u^\varepsilon(t)):=\int_{\T}
\vj(u^\varepsilon(\vx,t))d \vx
\equiv \int_{\T} \vj(u^\varepsilon_0(\vx))
 d \vx=\bvec{Q}(u^\varepsilon_0):=\bvec{Q}_0^\varepsilon \qquad t\ge0;
\end{equation}
and  the {\sl energy} defined as \cite{RN8,RN23,LinXin1999}
\begin{equation}\label{eq:conserve E}
  E(u^\varepsilon(t)):=\int_{\T}
 e(u^\varepsilon(\vx,t)) d\vx
\equiv \int_{\T}
e(u^\varepsilon_0(\vx))
d\vx=E(u^\varepsilon_0), \quad t\ge0.
\end{equation}
Here we adopt the notations,  for any complex-valued function $v(\vx):\ \T \to {\mathbb C}$, its corresponding current $\vj(v)$, energy density $e(v)$ and
Jacobian $J(v)$ are defined as
\begin{equation}\label{jvevJv}
\begin{aligned}
&\vj(v):=\im (\overline{v}\nabla v), \quad
e(v):=\frac{1}{2}|\nabla v|^2+\frac{1}{4\varepsilon^2}(1-|v|^2)^2,\quad J(v)=\frac{1}{2}\nabla\cdot(\J \vj(v))=\im (\partial_x\overline{v}\,\partial_y v),
\end{aligned}
\end{equation}
where $\overline {v}$ and $\im (v)$ denote the complex conjugate and imaginary part of the function $v$, respectively, and $\J$ is a symplectic matrix given as
\[\J=\left(\begin{array}{cc}0&1\\-1&0 \end{array}\right).\]

The NLSE \eqref{eq:maineq}, also known as the Gross-Pitaevskii equation (GPE) \cite{Serfaty2017,BaoBuZhang},  has been widely used as a phenomenological model for superfluidity, such as
liquid helium \cite{Pitaevskii1961,Kagan,LinXin1999} and Bose-Einstein condensation \cite{BaoBuZhang}. A key signature of superfluidity is the appearance of quantized vortices which are  particle-like or topological defects. Quantized vortices in two dimensions are those particle-like defects whose centers are zeros of the order parameter or wave function, possessing localized phase singularity with the topological charge (also known as winding number, index, or circulation) being quantized. They have been
widely observed in many different physical systems, such as liquid Helium, atomic gases, nonlinear optics and type-II
superconductors \cite{Pitaevskii1961,Bewley,BSX}. Their study remains
one of the most important and fundamental problems since they were predicted by
Lars Onsager in 1947 in connection with superfluid Helium.

Several analytical and numerical studies have dealt with quantized vortex states of the NLSE \eqref{eq:maineq} and their interactions as well as \blackmark{the reduced dynamical laws} of quantized vortex lattice when $\varepsilon\to0$ \cite{BaoTang,NEU1990385,Zhang2007}.
For results in the whole space ${\mathbb R}^2$ or on bounded domains with either Dirichlet or homogeneous Neumann boundary conditions, we refer to \cite{Bethuel2008,RN8,RN23,Lin1998,LinXin1999,Serfaty2017,RN2,RN27,RN5,Lin19992,Lin1999,RN100,RN101,BT2,BZZ} and references therein. Based on mathematical analysis and numerical simulation results \cite{Neu1990407,NEU1990385,Moronescu1995}, for a quantized vortex with winding number $m$, when $m=\pm1$, it is dynamically (or structurally) stable; and when $|m|>1$, it is unstable.

For the NLSE \eqref{eq:maineq} on \blackmark{the torus}, due to the periodic-type boundary condition,  there can exist several isolated and distinct quantized vortices in the initial data $u_0^\varepsilon:=u_0^\varepsilon(\vx)$, while the winding number of each quantized vortex is either $+1$ or $-1$ \cite{RN8}. It is well-known in the literature that the total number of quantized vortices in the initial data has to be an even integer and half of them with winding number $+1$ and the other half of them with winding number $-1$ \cite{RN8}.  We assume that in $u_0^\varepsilon$
there are $2N (N\in{\mathbb N})$ isolated and distinct quantized vortices whose centers are located
at $\va^{0,\varepsilon}_1,\va^{0,\varepsilon}_2,\ldots,\va^{0,\varepsilon}_{2N}\in \T$
 with winding number $d_1, d_2,\ldots,d_{2N} \in\{\pm1\}$, respectively.
 Without loss of generality, we assume
 \begin{equation}\label{d1d2d3}
 d_1=\ldots=d_N=1,\ d_{N+1}=\ldots=d_{2N}=-1, \quad
 \va^{0,\varepsilon}_j\neq \va^{0,\varepsilon}_k, \ 1\le j\ne k\le 2N.\qquad
 \end{equation}
Assume
\begin{equation}\label{vanmom}
\va^{0}_j:=\lim_{\varepsilon\to0} \va^{0,\varepsilon}_j, \quad 1\le j\le 2N.
\end{equation}
Then one has
\begin{equation}\label{eq:relation between Q and q}
\bvec{Q}_0:=\lim_{\varepsilon\to0} \bvec{Q}(u^\varepsilon_0) =2 \pi\J \sum_{j=1}^{2N}d_j\va^0_{j}=2 \pi\J \left[\sum_{j=1}^{N}\va^0_{j}-\sum_{j=1}^{N}\va^0_{N+j}\right].
\end{equation}

Under the assumption of the vanishing momentum of the initial data, i.e.
\begin{equation}\label{vanmomnew}
\va^{0}_j\neq \va^{0}_k, \quad  1\le j\ne k\le 2N, \qquad
\bvec{Q}_0=\lim_{\varepsilon\to0}\bvec{Q}_0^\varepsilon={\bf 0}\in {\mathbb R}^2,
\end{equation}
and
\begin{equation}\label{con:convergence of J varphi}
J(u^\varepsilon_0)\stackrel{\varepsilon\to 0^+}{\longrightarrow} \pi\sum_{j=1}^{2N} d_j \delta_{\va^0_{j}} \ \text{in}\ \M(\T):=[W^{1,\infty}(\T)]',
\end{equation}
with $\delta(\vx)$ the Dirac delta function and $\delta_{\va^0_{j}}:=\delta_{\va^0_{j}}(\vx)$ defined as $\delta_{\va^0_{j}}(\vx)=
\delta(\vx-\va^0_{j})$ for $1\le j\le 2N$,
Colliander and Jerrard \cite{RN8} established \blackmark{the reduced dynamical law} of
quantized vortex of the NLSE \eqref{eq:maineq} with \eqref{init} when $\varepsilon\to0$: for $t\ge0$, there exists the $j$-th vortex path, denoted as
$\va^{\varepsilon}_j(t)$ satisfying
$\va^{\varepsilon}_j(0)=\va^{0,\varepsilon}_j$,
 in the solution $u^\varepsilon(\vx,t)$
of the NLSE \eqref{eq:maineq} originated from $\va^{0,\varepsilon}_j$,  for $j=1,\ldots,2N$. Denote
\begin{equation}\label{path}
\begin{aligned}
&\va_j:=\va_j(t):=\lim_{\varepsilon\to0}\va^{\varepsilon}_j(t), \quad
j=1,\ldots,2N, \\
&\va:=\va(t)=(\va_1(t),\ldots,\va_{2N}(t))^T\in (\T)^{2N}_*,
\end{aligned}
\qquad t\ge0
\end{equation}
with
\begin{equation}
  (\T)^{2N}_*:=\left\{(\va_1,\cdots,\va_{2N})^T\in (\T)^{2N}|\va_k\ne \va_m\ \text{for any} \ 1\le k\ne m\le 2N \right\}.\end{equation}
Then when $\varepsilon\to0$, $\va(t)$ satisfies the following \blackmark{reduced dynamical law}:
\begin{equation} \label{oldODE}
\dot\va_j=-d_j\frac{1}{\pi}\J\nabla_{\va_j}W(\va)=2\J\sum_{1\le k\le 2N,\; k\ne j}d_k\, \nabla F(\va_j-\va_k),\qquad 1\le j\le 2N,
\end{equation}
with initial data
\begin{equation}\label{initODE}
\va_j(0)=\va_j^0, \qquad 1\le j\le 2N;
\end{equation}
where $W(\va)$ is the renormalized energy defined as
\begin{equation}\label{eq:define of W}
W(\va)=-\pi\sum_{1\le k\ne m\le 2N}d_kd_m\, F(\va_k-\va_m),
\end{equation}
with $F\in C^\infty_{loc}(\T\setminus\{\bvec{0}\})\cap W^{1,1}(\T)$ the solution of
\begin{equation}\label{eq:define of F}
\Delta F(\vx)=2\pi (\delta(\vx)-1),\ \vx\in \T, \quad {\rm with}\quad \int_{\T}F(\vx)d\vx=0.
\end{equation}

\begin{figure}[htp!]
\begin{center}
\begin{tabular}{cccc}
\includegraphics[height=3cm]{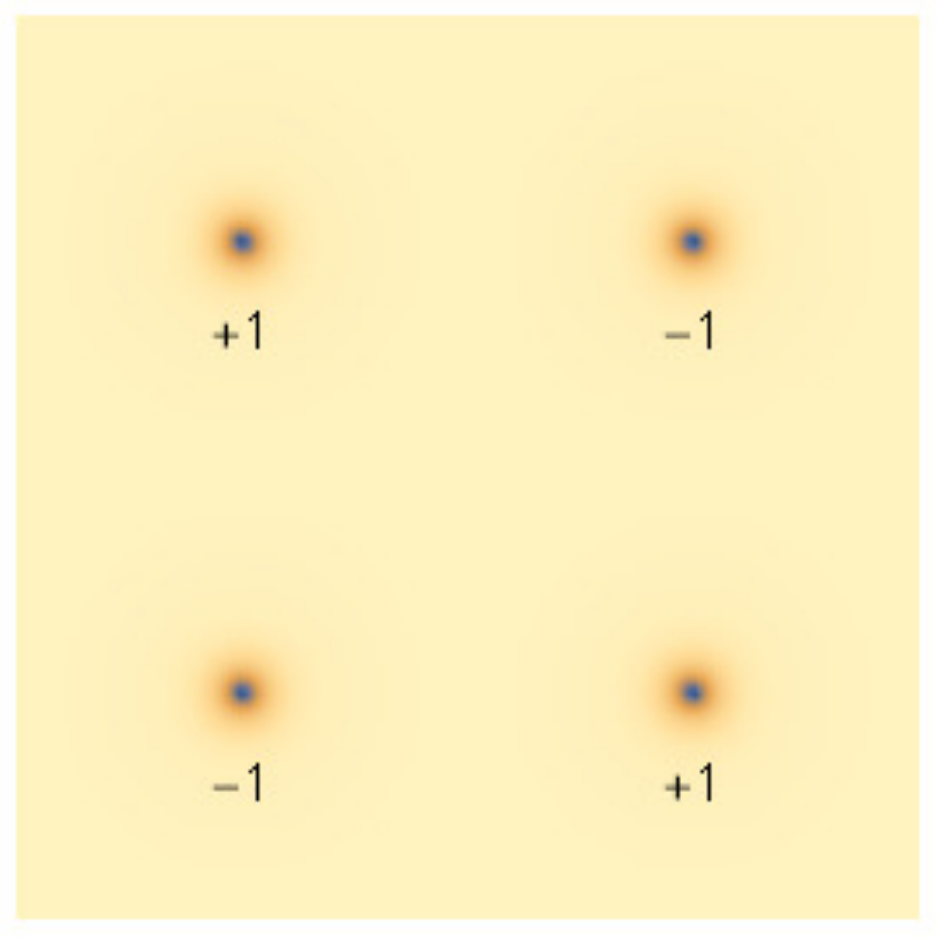}
&\includegraphics[width=3cm]{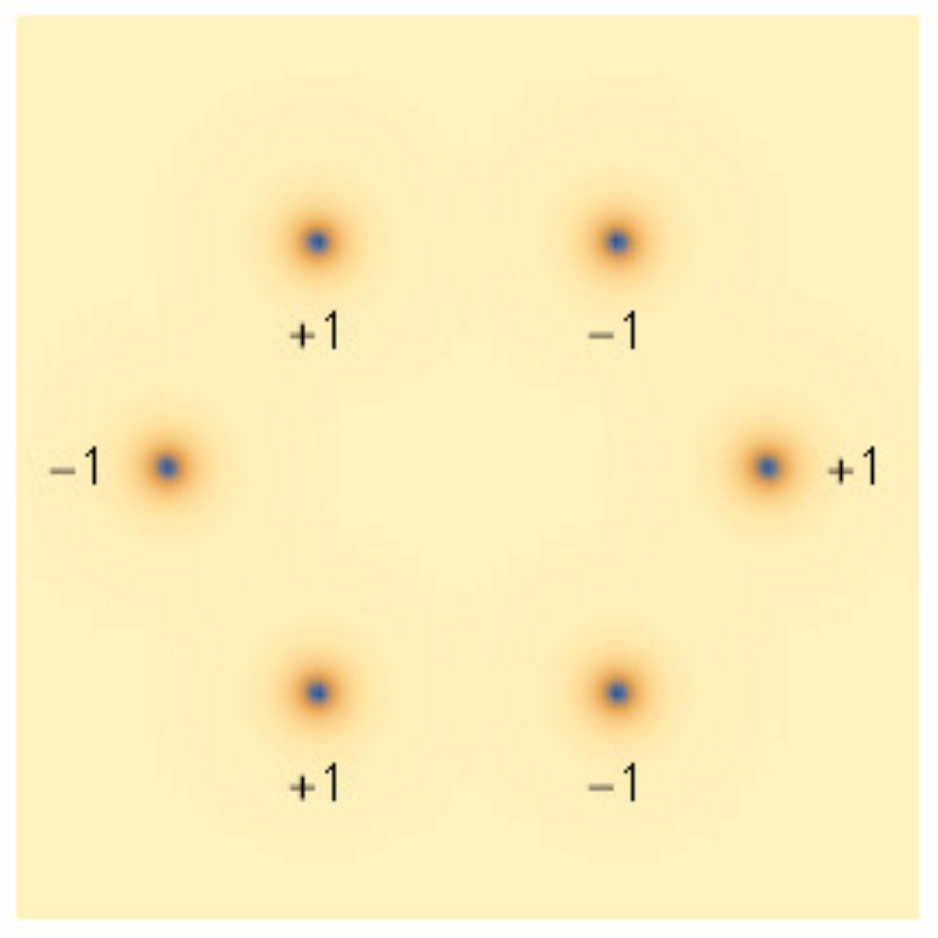}
&\includegraphics[width=3cm]{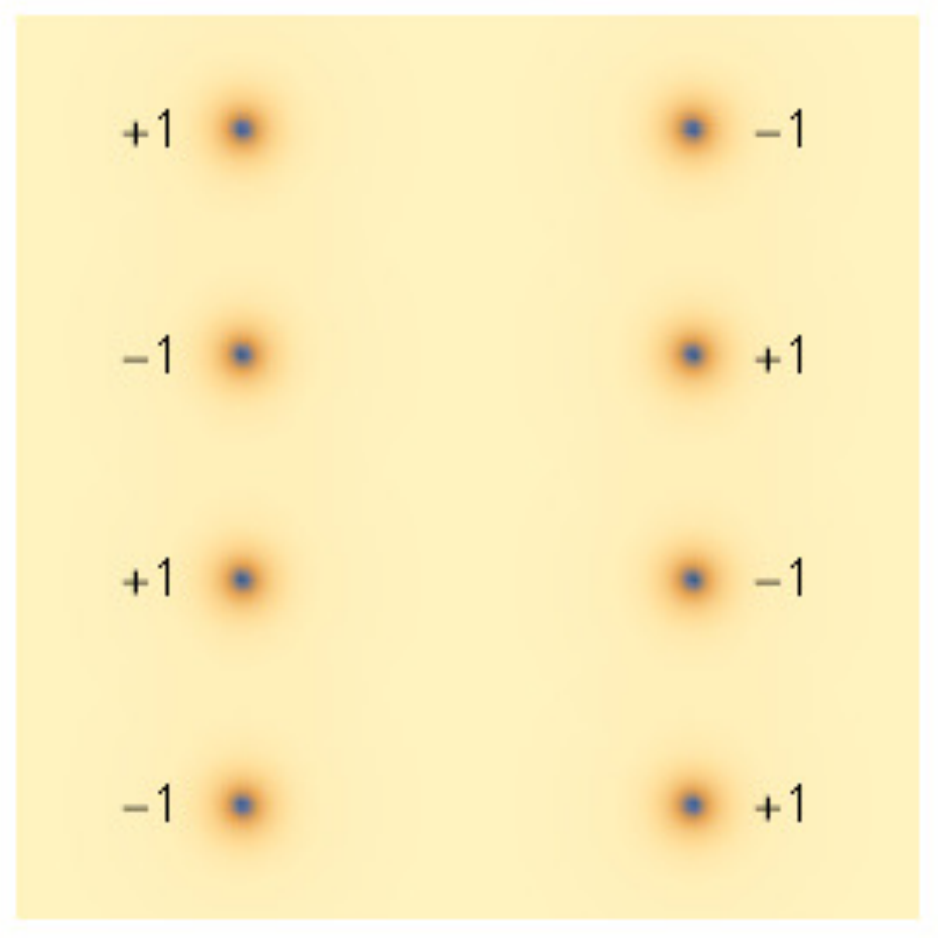}
&\includegraphics[height=3cm]{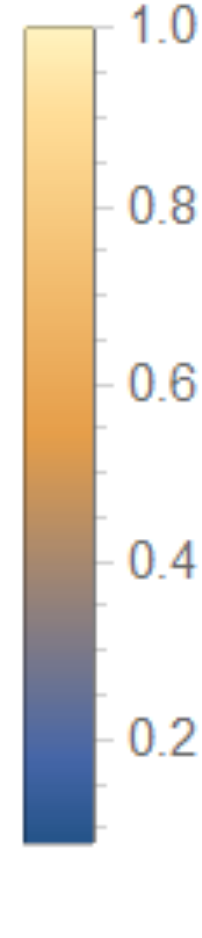}\\
\includegraphics[height=3cm]{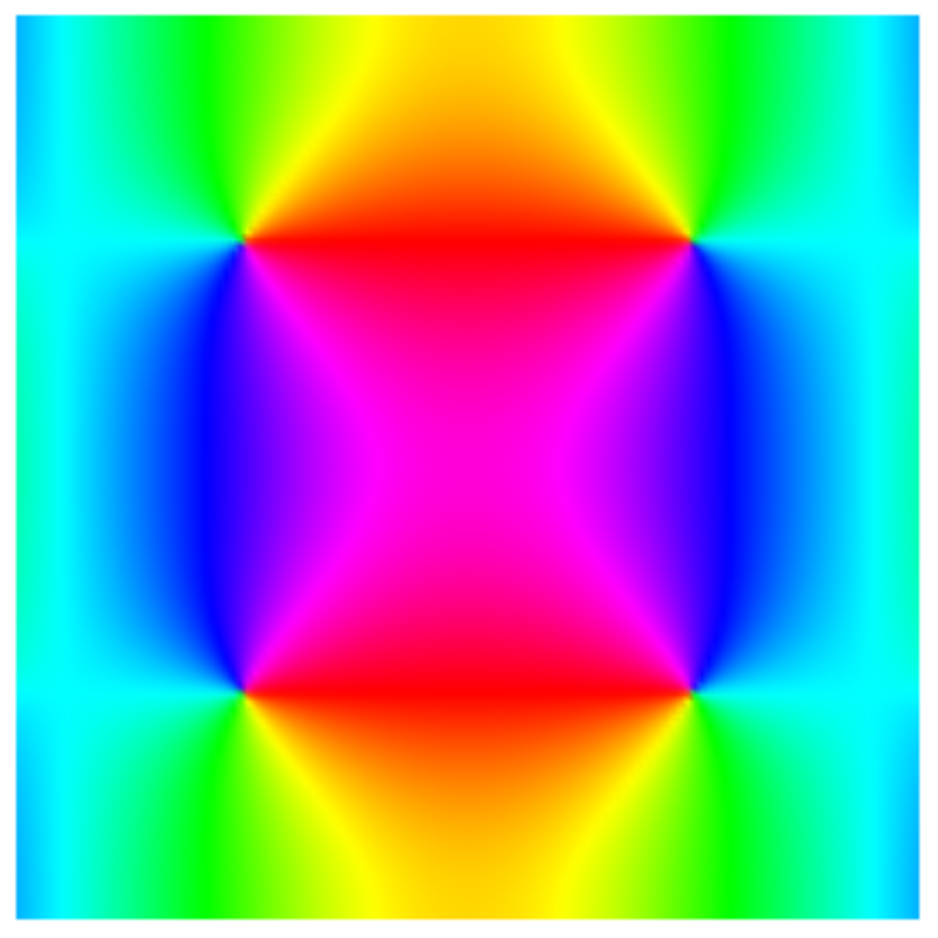}
&\includegraphics[width=3cm]{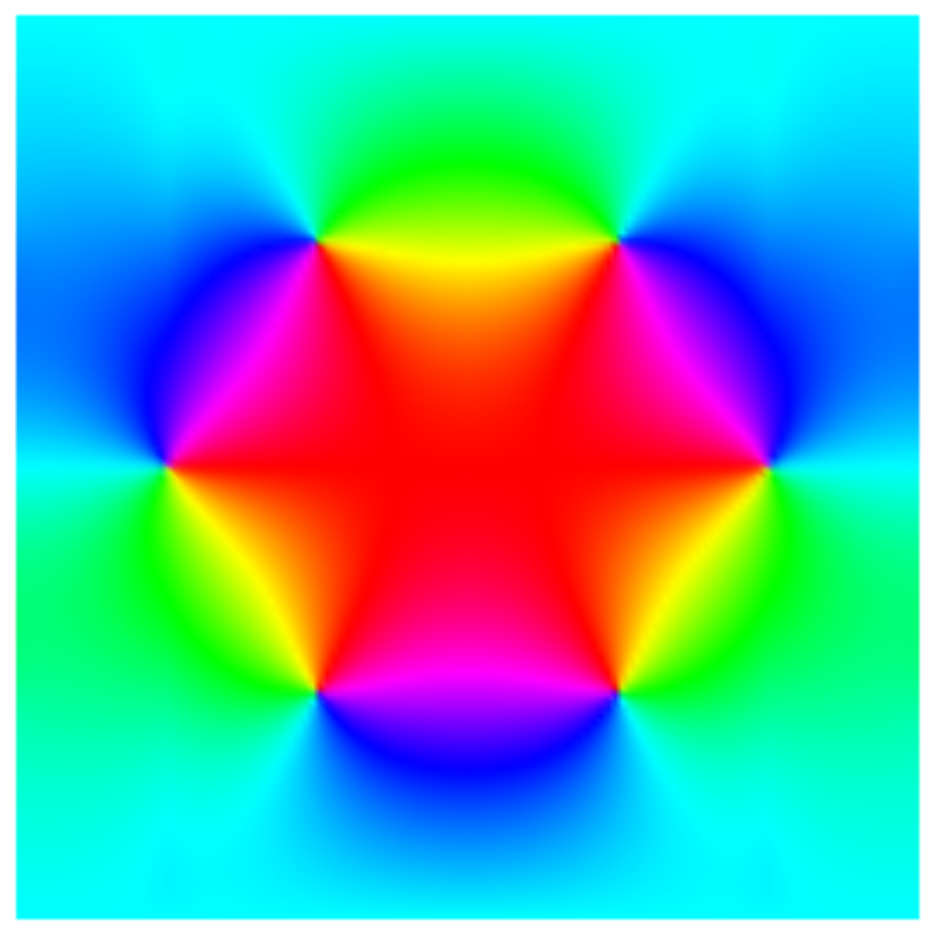}
&\includegraphics[width=3cm]{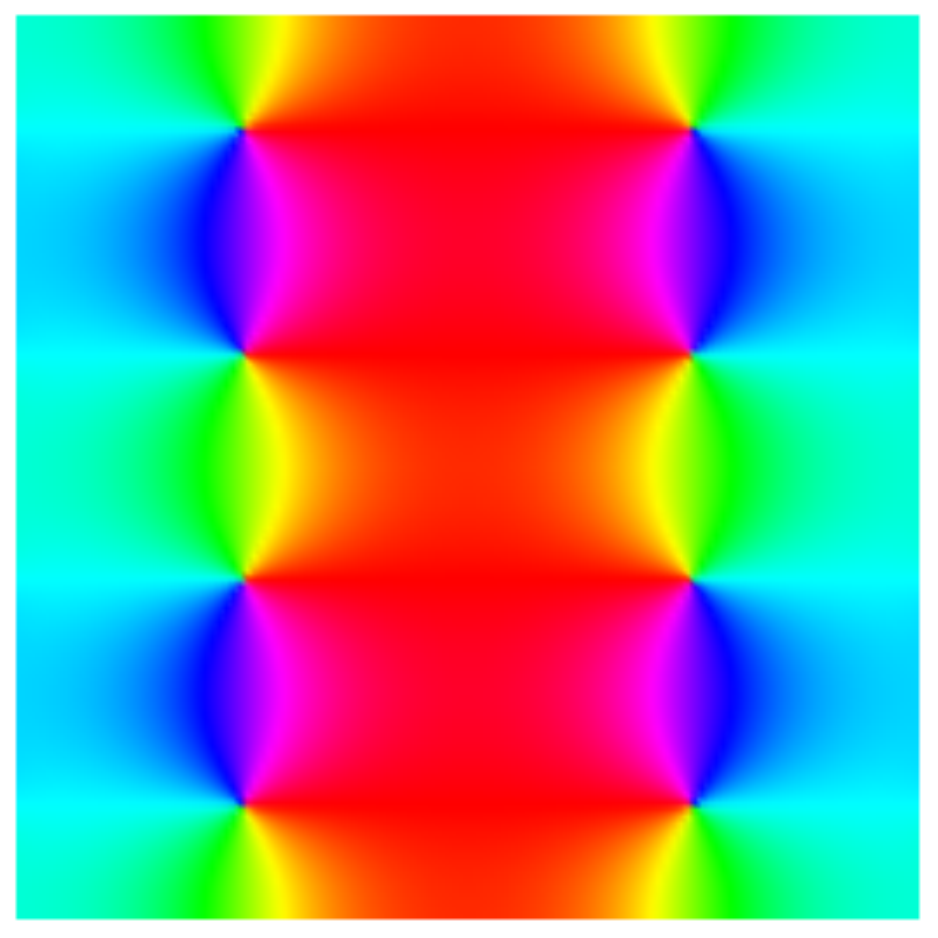}
&\includegraphics[height=3cm]{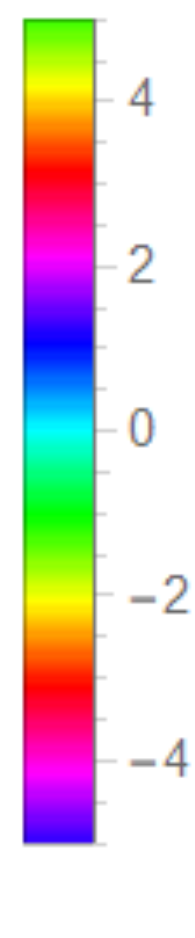}\\
\end{tabular}
\end{center}
\caption{Some typical vanishing momentum initial data $u_0^\varepsilon:=\sqrt{\rho_0^\varepsilon}e^{iS_0^\varepsilon}$:
Contour plots of  $|u^\varepsilon_0|$ with vortex center locations
with $+1$ and $-1$ as their corresponding winding numbers
(top row), and  contour plots of the phase function $S_0^\varepsilon$ (bottom row).
}\label{fig:examples of initial data vanishing}
\end{figure}

\begin{figure}[htp!]
\begin{center}
\begin{tabular}{cccc}
\includegraphics[height=3cm]{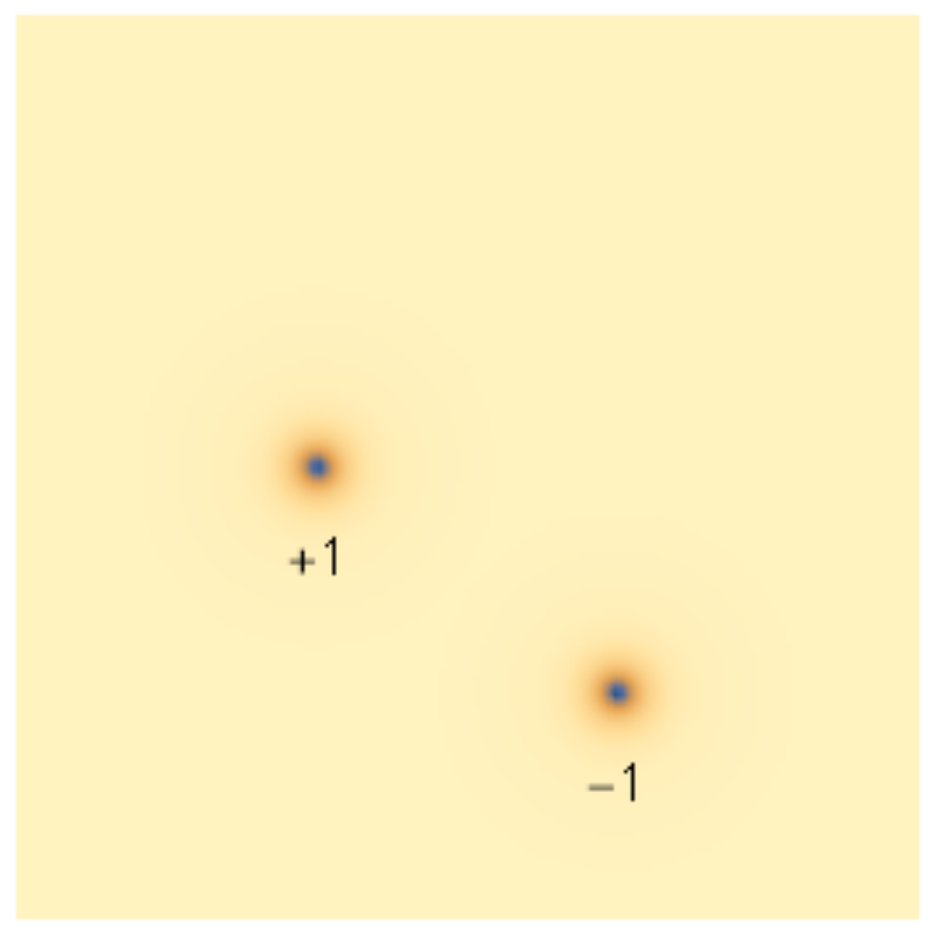}
&\includegraphics[width=3cm]{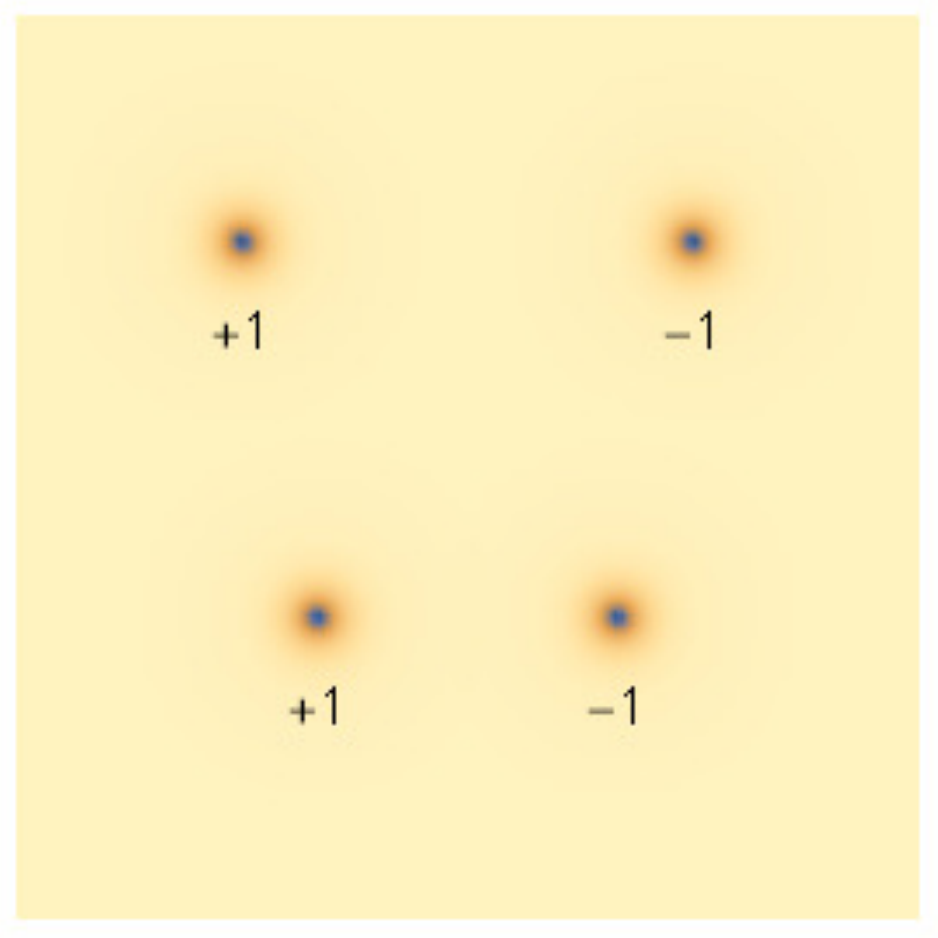}
&\includegraphics[width=3cm]{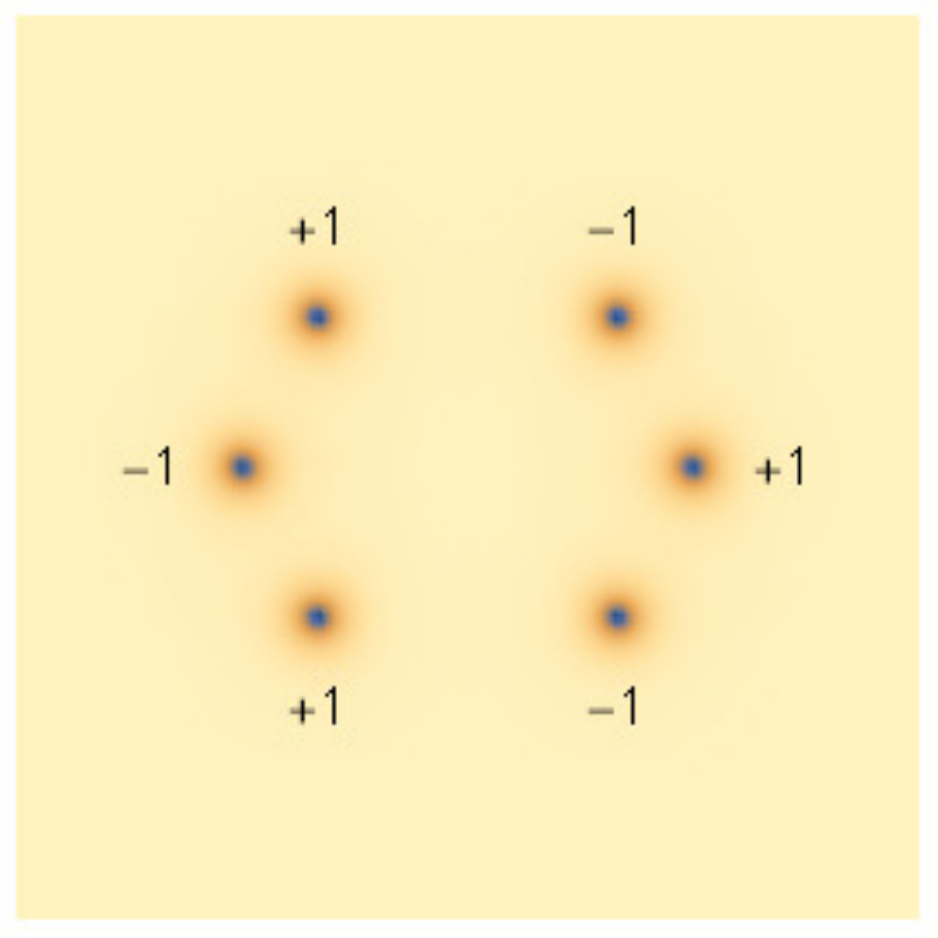}
&\includegraphics[height=3cm]{figures/normlegend.pdf}\\
\includegraphics[height=3cm]{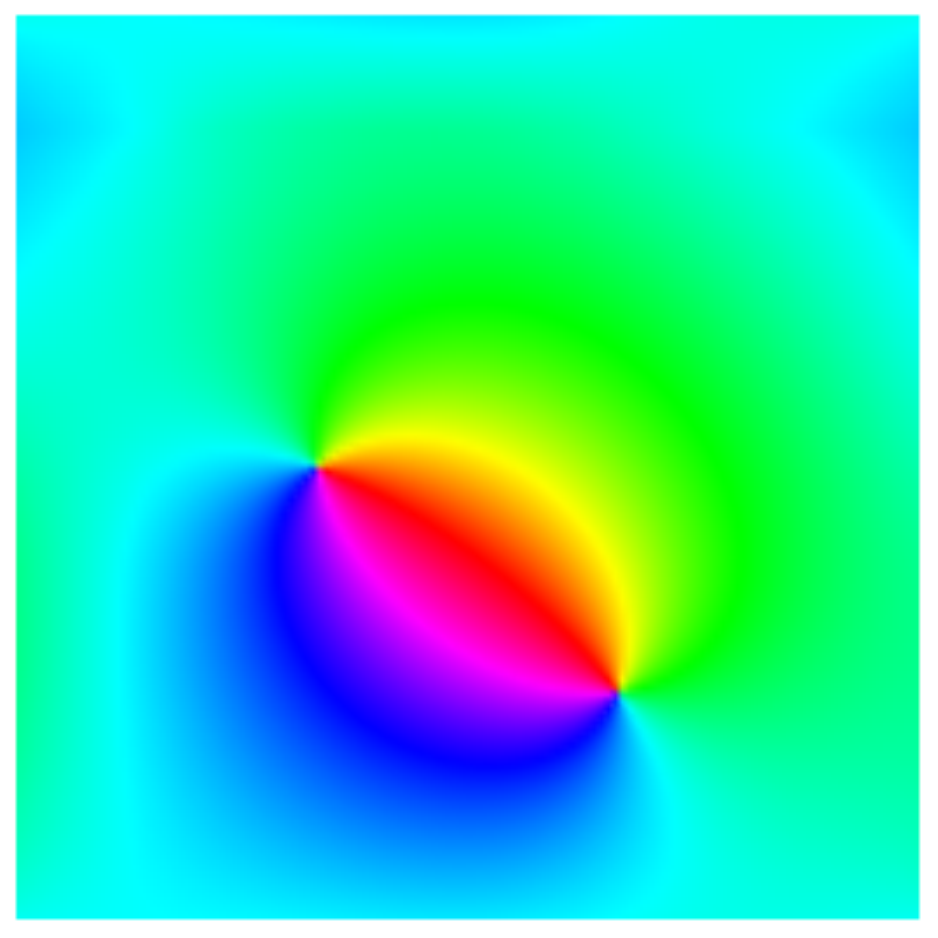}
&\includegraphics[width=3cm]{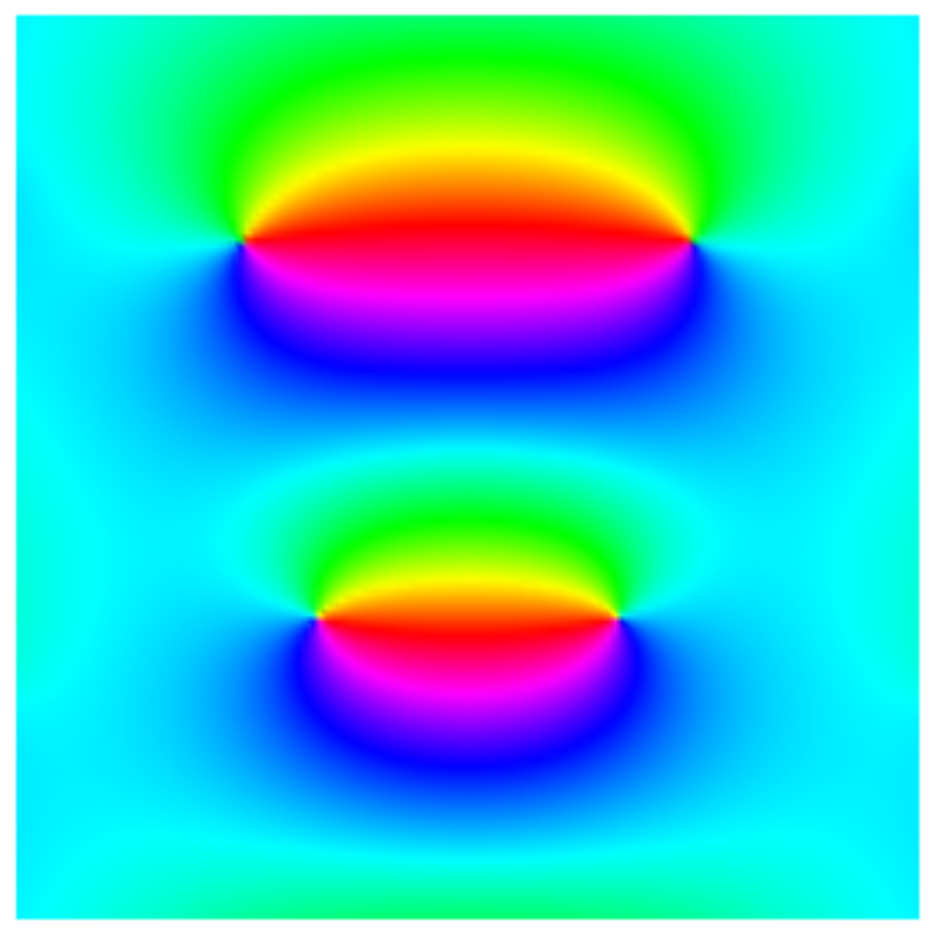}
&\includegraphics[width=3cm]{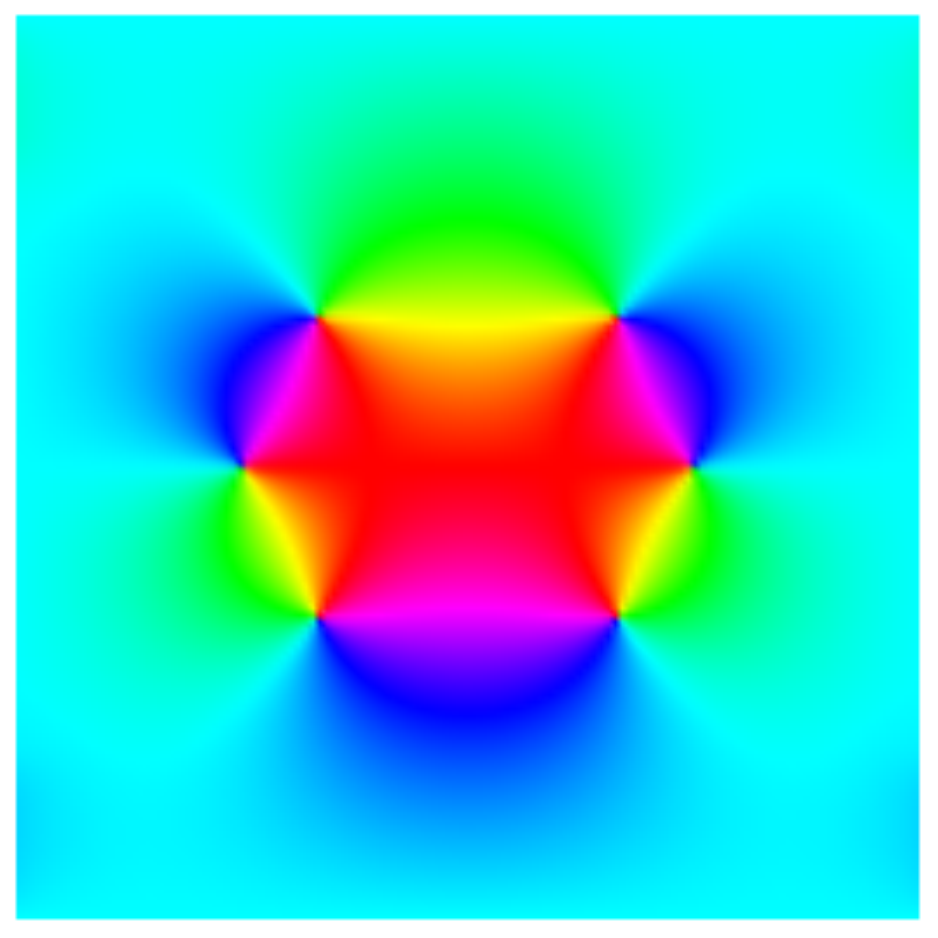}
&\includegraphics[height=3cm]{figures/phaselegend.pdf}\\
\end{tabular}
\end{center}
\caption{Some typical non-vanishing momentum initial data.
}\label{fig:examples of initial data nonvanishing}
\end{figure}

By \eqref{eq:relation between Q and q}, the vanishing momentum assumption used in \cite{RN8}, i.e. $\bvec{Q}_0={\bf 0}$, is equivalent to
\begin{equation}
\bvec{0}=\sum_{j=1}^{2N}d_j\va_j^0=\sum_{j=1}^{N}\va_j^0-\sum_{j=1}^{N}\va_{N+j}^0.
\end{equation}
Thus the vanishing momentum assumption is equivalent to
the assumption that
the positive mass center is the same as the negative mass center, i.e.
\begin{equation}\label{posnegcent}
\va_+^0=\va_-^0, \qquad {\rm with}\qquad \va_+^0:=\frac{1}{N}\sum_{j=1}^N\va_j^0, \quad \va_-^0:=\frac{1}{N}\sum_{j=N+1}^{2N}\va_j^0.
\end{equation}
From \eqref{vanmomnew} and \eqref{posnegcent}, one gets that $N\ge2$, i.e. $2N\ge4$, since $\va_+^0\ne\va_-^0$ when $N=1$.
In other words, \blackmark{the reduced dynamical law} for the NLSE \eqref{eq:maineq} obtained in Colliander and Jerrard \cite{RN8} only works for the case when the initial data
$u_0^\varepsilon$ admits $2N$ ($N\ge2$, i.e. $4$ or more)
vortices with their centers
satisfying proper symmetry requirement as stated in \eqref{posnegcent}.
Clearly, the requirement \eqref{posnegcent} on the initial data $u_0^\varepsilon$ excludes the non-vanishing momentum initial data, i.e.
$\bvec{Q}_0=\lim_{\varepsilon\to0}\bvec{Q}_0^\varepsilon\neq {\bf 0}$, which includes
two different cases: (i) $N=1$, i.e. two vortices; and (ii) $N\ge2$ and $\va_+^0\ne\va_-^0$. For the convenience of readers,
Figure \ref{fig:examples of initial data vanishing} shows some typical initial data with vanishing momentum, while Figure \ref{fig:examples of initial data nonvanishing} illustrates some typical initial data with non-vanishing momentum.

The main aim of this paper is to extend \blackmark{the reduced dynamical law} for quantized vortex dynamics of the NLSE \eqref{eq:maineq} with vanishing momentum initial data, i.e. $\bvec{Q}_0={\bf 0}$ in \eqref{vanmom}, to non-vanishing momentum initial data,
i.e. $\bvec{Q}_0\neq{\bf 0}$. A key ingredient is to construct a new canonical harmonic map to include the effect of the non-vanishing momentum into the dynamics.
To present our main result,  define
\begin{equation}\label{eq:def of I}
    I(\varepsilon)=\inf_{v\in H^1_g(B_1(\bvec{0}))}\int_{B_1(\bvec{0})}\left[\frac{1}{2}|\nabla v|^2+\frac{1}{4 \varepsilon^2}(|v|^2-1)^2\right]d\vx,
\end{equation}
where the function space $H^1_g(B_1(\bvec{0}))$ is defined as
\begin{align*}
  H_g^1(B_1(\bvec{0}))=\left\{v\in H^1(B_1(\bvec{0}))\left| v(\vx)=g(\vx)=\frac{x+iy}{|\vx|}\ \text{for}\  \vx\in \p B_1(\bvec{0})\right.\right\}.
\end{align*}
Define
\begin{equation}
  \gamma=\lim_{\varepsilon\to 0}\left(I(\varepsilon)-\pi\log\frac{1}{\varepsilon}\right).
\end{equation}
Introduce a  renormalized energy on \blackmark{the torus} as \cite{IgnatJerrard2021}
\begin{equation}\label{eq:define of WQ0}
\begin{aligned}
\W(\va)&=W(\va)+2\pi^2\left|\sum_{m=1}^{2N} d_m\va_m \right|^2\\&=-\pi\sum_{1\le k\ne m\le 2N}d_kd_m\, F(\va_k-\va_m) +2\pi^2\left|\sum_{m=1}^{2N} d_m\va_m \right|^2,\quad \va \in (\T)^{2N}_*,
\end{aligned}
\end{equation}
and  an $\varepsilon$-dependent renormalized energy as
\begin{equation}\label{eq:def of We}
\W^\varepsilon(\va):=2N\left(\pi\log\frac{1}{\varepsilon}+\gamma \right)+\W(\va),
\quad \va \in (\T)^{2N}_*.
\end{equation}

Our main result is stated as following:
\begin{Thm}[\blackmark{Reduced dynamical law} for NLSE]\label{thm:dynamics}
Assume the initial data $u_0^\varepsilon$ in \eqref{init} satisfies
\eqref{con:convergence of J varphi} and
\begin{equation}
 \label{Con:strong bd of energy}
\lim_{\varepsilon\to0} \bvec{Q}(u^\varepsilon_0)=\bvec{Q}_0:=2 \pi\J \left[\sum_{j=1}^{N}\va^0_{j}-\sum_{j=1}^{N}\va^0_{N+j}\right], \qquad
\limsup_{\varepsilon\to 0}\left[E(u^\varepsilon_0)-\W^\varepsilon(\val)\right]\le 0.
 \end{equation}
Then there exists a time $T>0$ and $2N$ Lipschitz paths given in \eqref{path}
such that the solution $u^\varepsilon$ of the NLSE \eqref{eq:maineq} with
\eqref{init} satisfies
\begin{equation}\label{eq:con of J}
J(u^\varepsilon(t))\stackrel{\varepsilon\to 0^+}{\longrightarrow} \pi\sum_{j=1}^{2N} d_j \delta_{\va_{j}(t)} \ \text{in}\ \M(\T),
\end{equation}
where $\va_j(t)$ ($1\le j\le 2N$) satisfy the following \blackmark{reduced dynamical law}:
\begin{equation}\label{newODE}
  \dot{\va}_j=
-d_j\frac{1}{\pi}\J\nabla_{\va_j}W(\va)-2\bvec{Q}_0=2\J\sum_{1\le k\le 2N,\; k\ne j}d_k\, \nabla F(\va_j-\va_k)-2\bvec{Q}_0,\quad  t>0,
\end{equation}
with the initial data \eqref{initODE}.
\end{Thm}

We remark here that: (i) when $\bvec{Q}_0={\bf 0}$, \eqref{newODE} collapses to \eqref{oldODE}; and (ii) \eqref{newODE} is also equivalent to the following
ordinary differential equations (ODEs)
\begin{equation}\label{eq:dynamics}
\dot{\va}_j=-d_j\frac{1}{\pi}\J\nabla_{\va_j}\W(\va)=
2\J\sum_{1\le k\le 2N,\; k\ne j}d_k\, \nabla F(\va_j-\va_k)-4\pi\vq(\va),\quad  t>0,
\end{equation}
where $\vq(\va)$ is defined as
\begin{equation}\label{vqva}
\vq(\va):= \J\sum_{m=1}^{2N}d_m\va_m, \qquad \va\in (\T)^{2N}_*.
\end{equation}
Since $\vq(\va)$ is a first integral of \eqref{eq:dynamics}, and if $\va$ is the solution
of \eqref{eq:dynamics} with the initial data \eqref{initODE}, we have
\begin{equation}
\vq(\va(t))\equiv \vq(\va(0))=\J\sum_{m=1}^{2N}d_m\va_m^0=\frac{1}{2\pi}\bvec{Q}_0,
\qquad t\ge0.
\end{equation}

The paper is organized as follows. In Section \ref{sec:pre}, we introduce a new
canonical harmonic map on \blackmark{the torus} corresponding to $2N$ vortex centers with non-vanishing initial momentum $\bvec{Q}_0\ne {\bf 0}$ and show its properties.
In Section \ref{sec:exist}, we prove the local existence of $2N$ vortex paths in the solution of the NLSE and the convergence of the corresponding current. In Section \ref{sec:proof main}, we establish \blackmark{the reduced dynamical law} \eqref{eq:dynamics}.
In Section \ref{sec:pro of rdl}, we present some first integrals and give \blackmark{several} analytical solutions of \blackmark{the reduced dynamical law} with initial data with symmetry. Finally, some concluding remarks are drawn in Section \ref{sec:conclusion}.

\section{Canonical harmonic maps}\label{sec:pre}
To prove \blackmark{the reduced dynamical law} \eqref{newODE} in Theorem \ref{thm:dynamics}, similar to the proof in \cite{RN8} for the case of vanishing initial momentum, i.e. $\bvec{Q}_0= {\bf 0}$, we introduce a new canonical harmonic map on \blackmark{the torus} which works for both non-vanishing initial momentum, i.e.
$\bvec{Q}_0\ne {\bf 0}$, and vanishing initial momentum, i.e. $\bvec{Q}_0= {\bf 0}$,
by adopting a phase shift depending on $\bvec{Q}_0$.

\subsection{A canonical harmonic map for a vortex dipole}\label{sec:construction}\label{sec:def of H v2}
For simplicity of notations, we first consider the most simple case of a vortex dipole, i.e. $N=1$, and assume the two vortex centers are located at $\va_1=(x_1,y_1)^T\ne\va_2=(x_2,y_2)^T\in\T$
with winding number $d_1=1$ and $d_2=-1$, respectively. Denote
\begin{equation}\label{Omqa}
\Omega=[0,1]^2, \qquad
\va=(\va_1,\va_2)^T\in(\T)^2, \qquad \vq_2(\va):=(q_1,q_2)^T=\J(\va_1-\va_2).
\end{equation}
Using the known results on harmonic maps in Section I.3 in \cite{Bethuel2017}, there exists a complex-valued function
$\widetilde{H}:=\widetilde{H}(\vx;\va)\in C_{\rm loc}^\infty(\Omega\setminus\{\va_1,\va_2\})\cap W^{1,1}(\Omega)$ \blackmark{that satisfies}
\begin{equation}\label{vjtildeH}
  \vj(\blackmark{\widetilde{H}(\vx;\va)})=-\J\nabla(F(\vx-\va_1)-F(\vx-\va_2)), \quad\
   |\widetilde{H}(\vx;\va)|=1,
   \quad \vx\in \Omega\setminus\{\va_1,\va_2\}.
\end{equation}
Unfortunately, $\widetilde{H}\not\in C^\infty_{\rm loc}(\T\setminus\{\va_1,\va_2\})\cap W^{1,1}(\T)$ \blackmark{due to} $\va_1-\va_2\ne{\bf 0}$ \cite{Bethuel2017,RN8}.

\begin{Lem}[A canonical harmonic map for a vortex dipole]\label{lem:def of H v2}
Define a complex-valued function as
\begin{equation}\label{H2def}
  H_2:=H_2(\vx;\va)=\blackmark{\widetilde{H}(\vx;\va)}e^{2\pi i\vq_2(\va)\cdot\vx}, \quad \vx\in \T\setminus\{\va_1,\va_2\}.
\end{equation}
Then $H_2\in C^\infty_{\rm loc}(\T\setminus\{\va_1,\va_2\})\cap W^{1,1}(\T)$ satisfies
$|H_2(\vx;\va)|=1$ and
\begin{equation}\label{eq:def of H 2v}
  \vj(\blackmark{H_2(\vx;\va)})=-\J\nabla(F(\vx-\va_1)-F(\vx-\va_2))+2\pi\vq_2(\va),
   \quad \vx\in \T\setminus\{\va_1,\va_2\}.
\end{equation}
\end{Lem}

\begin{proof} From $\widetilde{H}:=\widetilde{H}(\vx;\va)\in C_{\rm loc}^\infty(\Omega\setminus\{\va_1,\va_2\})\cap W^{1,1}(\Omega)$ and $|\widetilde{H}(\vx;\va)|=1$, noting
\eqref{H2def}, we know that $H_2(\vx;\va)\in C_{\rm loc}^\infty(\Omega\setminus\{\va_1,\va_2\})\cap W^{1,1}(\Omega)$ and
$|H_2(\vx;\va)|=1$. Combining \eqref{jvevJv}, \eqref{H2def}, \eqref{vjtildeH}  and \eqref{Omqa}, we obtain
\begin{equation}\label{eq:splitting of j(HH)}
\begin{aligned}
\vj(\blackmark{H_2(\vx;\va)})=&\im (\blackmark{\overline{H_2(\vx;\va)}}\nabla \blackmark{H_2(\vx;\va)})\\=&\im\left(\blackmark{\overline{\widetilde{H}(\vx;\va)}}e^{-2\pi i\vq_2(\va)\cdot\vx}e^{2\pi i\vq_2(\va)\cdot\vx}(\nabla\blackmark{\widetilde{H}(\vx;\va)}+2\pi i\blackmark{\widetilde{H}(\vx;\va)}\vq_2(\va) )\right)\\=&\vj(\blackmark{\widetilde{H}(\vx;\va)})+2\pi\vq_2(\va)=-\J\nabla(F(\vx-\va_1)-F(\vx-\va_2))+2\pi\vq_2(\va),
\end{aligned}
\end{equation}
which immediately implies \eqref{eq:def of H 2v}.

In order to show the periodicity of $H_2$, noting $|H_2(\vx;\va)|=1$, we have
\begin{equation}\label{H2Phi}
  H_2(\vx;\va)=e^{i\Theta(\vx)}, \qquad \vx\in \T\setminus\{\va_1,\va_2\},
\end{equation}
where $\Theta(\vx)$ is the phase function. Combining \eqref{jvevJv}, \eqref{H2Phi},
\eqref{eq:def of H 2v}  and
\eqref{Omqa}, we get
\begin{equation}\label{eq:derivative of Theta}
\begin{aligned}
\vj(\blackmark{H_2(\vx;\va)})&=\im (\blackmark{\overline{H_2(\vx;\va)}}\nabla \blackmark{H_2(\vx;\va)})=\im (e^{-i\Theta(\vx)}e^{i\Theta(\vx)}i\nabla\Theta(\vx) )\\
&=\nabla \Theta(\vx)=-\J\nabla(F(\vx-\va_1)-F(\vx-\va_2))+2\pi\vq_2(\va).
\end{aligned}
\end{equation}
From \eqref{eq:derivative of Theta}, \eqref{eq:define of F} and \eqref{Omqa},
we get $\nabla \Theta\in C^\infty_{\rm loc}(\T\setminus\{\va_1,\va_2\})\cap L^1(\T)$.

Without loss of generality, we can assume $x_1< x_2,y_1\le y_2$.
For the case $y_1<y_2$, from \eqref{eq:derivative of Theta} and \blackmark{noting that $F$} is a function defined on \blackmark{the torus} in \eqref{eq:define of F},  we have
\begin{align*}
\Theta(1,y)-\Theta(0,y)&=\int_0^1 \p_{x}\Theta\left(x,y \right)dx\\
&=-\int_0^1 \left[\p_y F\left(x-x_1,y-y_1\right)-\p_y F\left(x-x_2,y-y_2\right)\right]d x+2\pi q_1\\
&=-\int_0^1 \left[\p_y F\left(x,y-y_1\right)-\p_y F\left(x,y-y_2\right)\right]d x+2\pi q_1\\
&=-\int_{\p\Gamma}\frac{\p F}{\p\vn}ds+2\pi q_1=-\int_{\Gamma}\Delta Fd\vx+2\pi(y_1-y_2)\\
&=-\int_{\Gamma}2\pi(\delta(\vx)-1)d\vx+2\pi(y_1-y_2)=\left\{\begin{array}{cc}
  0,&y\not \in(y_1,y_2),\\-2\pi,&y\in(y_1,y_2) ,\end{array} \right.
\end{align*}
where $\Gamma=[0,1]\times[y_1,y_2]$ and $\vn$ is the unit outward normal vector. For the case $y_1=y_2$, we have
\begin{align*}
\Theta(1,y)-\Theta(0,y)=&\int_0^1\p_x\Theta(x,y)dx\\
=&-\int_0^1[\p_yF(x-x_1,y-y_1)-\p_yF(x-x_2,y-y_2)]dx+2\pi q_1\\
=&-\int_0^1[\p_yF(x,y-y_1)-\p_yF(x,y-y_2)]dx+2\pi(y_1-y_2)=0.
\end{align*}
Combining the above two equalities, we get
\begin{equation}\label{eq:PBCy}
    \Theta(1,y)-\Theta(0,y)=\int_0^1\p_{x}\Theta\left(x,y \right)dx=\left\{\begin{array}{ll}
    0,&y\notin(y_1,y_2),\\
    -2\pi,&y\in(y_1,y_2).
    \end{array}\right.
\end{equation}
Plugging \eqref{eq:PBCy} into \eqref{H2Phi}, we get
\begin{equation}\label{H21y}
\blackmark{H_2((1,y)^T;\va)}=e^{i\Theta(1,y)}=e^{i\Theta(0,y)}=\blackmark{H_2((0,y)^T;\va)}, \qquad 0\le y\le1.
\end{equation}
Similarly, we can prove
\begin{equation}\label{eq:PBCx}
    \Theta(x,1)-\Theta(x,0)=\int_0^1\p_{y}\Theta\left(x,y \right)dy=\left\{\begin{array}{ll}
    0,&x\notin(x_1,x_2),\\
    2\pi,&x\in(x_1,x_2),
    \end{array}\right.
\end{equation}
and
\begin{equation}\label{H2x1}
\blackmark{H_2((x,1)^T;\va)}=e^{i\Theta(x,1)}=e^{i\Theta(x,0)}=\blackmark{H_2((x,0)^T;\va)}, \qquad 0\le x\le1.
\end{equation}
Taking \blackmark{the gradient of} \eqref{H2Phi}, we have
\begin{equation}\label{gradH2}
\nabla H_2(\vx;\va)=i\, H_2(\vx;\va)\nabla \Theta(\vx).
\end{equation}
Combining \eqref{H2Phi}, \eqref{gradH2}, \eqref{H2x1} and \eqref{H21y}, and noting
$\nabla \Theta\in C^\infty_{\rm loc}(\T\setminus\{\va_1,\va_2\})\cap L^1(\T)$, we obtain immediately $H_2\in C^\infty_{\rm loc}(\T\setminus\{\va_1,\va_2\})\cap W^{1,1}(\T)$.
\end{proof}

\subsection[A canonical harmonic map for N vortex dipoles and its properties]{A canonical harmonic map for \texorpdfstring{$N$}{}  vortex dipoles and its properties}

For $2N$ vortex centers $\va_1,\cdots,\va_{2N}\in\T$ and $\va=(\va_1,\cdots,\va_{2N})^T\in(\T)^{2N}_*$, we can divide $\va$ into $N$ vortex dipoles: $(\va_1,\va_{N+1})^T,\cdots,(\va_N,\va_{2N})^T$. Then, we define the canonical harmonic map as
\begin{equation}\label{eq:def of H2}
  H:=H(\vx;\va)=\prod_{j=1}^NH_2(\vx;(\va_j,\va_{N+j})^T).
\end{equation}
Similar to \eqref{eq:splitting of j(HH)}, we have
\begin{align}\label{eq:def of H}
  \blackmark{\vj(H):=}&\vj(H(\vx,\va))=\sum_{j=1}^N\vj(H_2(\vx;(\va_j,\va_{N+j})^T))\nn
  =&\sum_{j=1}^{N}\left[-\J\nabla(F(\vx-\va_j)-F(\vx-\va_{N+j}))+
  2\pi\vq_2(\left(\va_j,\va_{N+j}\right)^T)\right] \nn=&-\J\sum_{j=1}^{2N}d_j\nabla F(\vx-\va_j)+2\pi\vq(\va).
\end{align}

As shown in Lemma \ref{lem:def of H v2}, we have $H_2(\vx;(\va_j,\va_{N+j})^T)\in C^\infty_{\rm loc}(\T\setminus\{\va_j,\va_{N+j}\})\cap W^{1,1}(\T)$ and \\$|H_2(\vx,(\va_j,\va_{N+j})^T)|=1$, which implies 
\begin{equation}\label{eq:regularity of H}
  H\in C^\infty_{\rm loc}(\blackmark{\T_*(\va)})\cap W^{1,1}(\T),\qquad |H(\vx;\va)|=1,
\end{equation}
with
\begin{equation}\label{eq:def of T_*a}
  \T_*(\va)=\T\setminus\{\va_1,\cdots,$ $\va_{2N}\}.
\end{equation}

We then introduce some notations. For $z,Z\in\mC$, we define
\[
  \la z,Z\ra:=\frac{1}{2}(\overline{z}Z+z\overline{Z}).
\]
And for two complex vectors $\bvec{z}=(z_1,z_2)^T,\bvec{Z}=(Z_1,Z_2)^T\in \mC^2,$ we define
\[
  \la\bvec{z},\bvec{Z}\ra:=\la z_1,Z_1\ra+\la z_2,Z_2\ra.
\]
In particular, if $\bvec{z},\bvec{Z}\in\R^2\subset\mC^2$,
\[
  \la\bvec{z},\bvec{Z}\ra=\bvec{z}\cdot\bvec{Z}.
\]
 We  denote
 \begin{equation}\label{eq:def of Phi}
\Psi(\vx):=\sum_{j=1}^{2N}d_jF(\vx-\va_j),
 \end{equation}
 and for $0<\rho\ll 1$
 \begin{equation}\label{eq:def of T_r}
    \T_\rho(\va):=\T\setminus\cup_{j}B_\rho(\va_j).
\end{equation}
In the following, we use $\Hess(\eta)$ to denote the Hessian matrix of a function $\eta$.

We then derive some simple properties of $H$ defined in \eqref{eq:def of H}:
\begin{Lem}For the  $H=H(\vx;\va)$ given in \eqref{eq:def of H2}, we have
\begin{align}
\Div \vj(H)=0,\qquad J(H)=\pi\sum_{j=1}^{2N}d_j\delta_{\va_j},\qquad \int_{\T}\vj(H)d\vx=&2\pi\bvec{q}(\va),\quad \va \in (\T)^{2N}_*, \label{eq:div of jH}
\end{align}
and for $0<\rho\ll1 $
\begin{align}\label{eq:est of int of |jH|^2}
 \int_{\T_\rho(\va)}e(H)d \vx= 2N \pi \log \frac{1}{\rho} +\W(\va)+O\left(\rho^2 \right),\quad \va \in (\T)^{2N}_*.
\end{align}
\end{Lem}

\begin{proof}
Via direct calculation of $\nabla\cdot\vj(H)$ and $J(H)$, noting \eqref{jvevJv},  \eqref{eq:def of H} and that $F\in W^{1,1}(\T)$ satisfies \eqref{eq:define of F}, we have
\begin{equation}
  \nabla\cdot \vj(\blackmark{H(\vx;\va)})=-\nabla\cdot\left(\sum_{j=1}^{2N}d_j\J\nabla F(\vx-\va_j)\right)=0,
\end{equation}
\begin{equation}
\begin{aligned}
  J(\blackmark{H(\vx;\va)})=&\frac{1}{2}\nabla\cdot(\J\vj(\blackmark{H(\vx;\va)}))=-\frac{1}{2}\nabla\cdot\left(\J\sum_{j=1}^{2N}d_j\J\nabla F(\vx-\va_j)\right)\\=&\frac{1}{2}\sum_{j=1}^{2N}d_j\Delta F(\vx-\va_j)=\pi\sum_{j=1}^{2N}d_j\delta_{\va_j}(\vx).
  \end{aligned}
\end{equation}
Integrating $\vj(H)$ over $\T$ and noting \eqref{eq:def of H}, we obtain
  \begin{align}
  \int_{\T}\vj(H)d\vx=&-\J\sum_{j=1}^{2N}d_j\int_{\T}\nabla F(\vx-\va_j)d\vx+2\pi\int_{\T}\vq(\va) d\vx\nn
  =&-\J\sum_{j=1}^{2N}d_j\int_{\T}\nabla F(\vx)d\vx+2\pi\vq(\va)=2\pi\vq(\va).
\end{align}
Combining the above three equalities, we obtain \eqref{eq:div of jH}.

Noting $|H(\vx;\va)|=1$, we can assume $\blackmark{H(\vx;\va)}=e^{i\theta(\vx)}$. Similar to \eqref{eq:derivative of Theta} and \eqref{gradH2}, we have
\begin{equation}\label{eq:j(H)=DH}
   \vj(H)=\nabla \theta ,\quad \nabla H=iH\nabla\theta=iH\vj(H).
 \end{equation}
  Then plugging $|H|=1$, \eqref{eq:j(H)=DH}, \eqref{eq:def of Phi} and \eqref{eq:def of H} into the definition of $e(H)$ in \eqref{jvevJv}, and integrating $e(H)$ over $\T_\rho(\va)$, we have
 \begin{align}\label{eq:est of int of |jH|^22}
 \int_{\T_\rho(\va)}e(H)d\vx=&\frac{1}{2}\int_{\T_\rho(\va)}|\nabla H|^2d \vx=\frac{1}{2}\int_{\T_\rho(\va)}|\vj(H)|^2d \vx=\frac{1}{2}\int_{\T_\rho(\va)}|-\J\nabla\Psi+2\pi\bvec{q}(\va)|^2d \vx\nn
 =& \frac{1}{2}\int_{\T_\rho(\va)}|\nabla\Psi|^2d\vx-2\pi\int_{\T_\rho(\va)}\vq(\va)\cdot(\J\nabla\Psi)d\vx+2\pi^2|\vq(\va)|^2+O(\rho^2).
\end{align}
Similar to Lemma 12 in \cite{RN23}, the first term on the right-hand side of \eqref{eq:est of int of |jH|^22} can be estimated by
\begin{equation}\label{eq:first term est}
  \frac{1}{2}\int_{\T_\rho(\va)}|\nabla\Psi|^2d\vx = W(\va)+2N\pi\log\frac{1}{\rho}+O(\rho^2).
\end{equation}
For the second term on the right-hand side of \eqref{eq:est of int of |jH|^22}, we  define
\begin{equation}\label{eq:psij}
\Psi_j(\vx)=\Psi(\vx)-d_j\log|\vx-\va_j|.
\end{equation}
Then  $\Psi_j\in C^1(B_\rho(\va_j))$. Noting that $\J\vq(\va)$ is constant with respect to $\vx$ and substituting \eqref{eq:psij} into the second term on the right-hand side of \eqref{eq:est of int of |jH|^22}, we have
\begin{align}\label{eq:second term est}
\int_{\T_\rho(\va)}\vq(\va)\cdot(\J\nabla\Psi)d\vx=&-\int_{\T_\rho(\va)}(\J\vq(\va))\cdot\nabla\Psi d\vx
=\sum_{j=1}^{2N}\int_{\p B_\rho(\va_j)}\bvec{n}\cdot(\J\vq(\va))\Psi ds\nn
=&\sum_{j=1}^{2N}\int_{\p B_\rho(\va_j)}(d_j\log|\vx-\va_j|+\Psi_j(\vx))(\J\vq(\va))\cdot\bvec{n}ds\nn
=&\sum_{j=1}^{2N}\left(d_j\log\rho\int_{\p B_\rho(\va_j)}(\J\vq(\va))\cdot \bvec{n}ds+\int_{B_\rho(\va_j)}\nabla \Psi_j(\vx)\cdot(\J\vq(\va))d\vx\right)\nn
=&O(\rho^2).
\end{align}
Substituting \eqref{eq:first term est}, \eqref{eq:second term est}, \eqref{eq:define of WQ0} and \eqref{vqva}  into \eqref{eq:est of int of |jH|^22}, we immediately obtain \eqref{eq:est of int of |jH|^2}.
\end{proof}

\begin{Lem}\label{lem:production of jH and eta}
Suppose that $\eta\in C^2(\T)$ is linear in a neighborhood of $\va_{j}$ and
$\supp (\eta)\cap \{ \va_1,\cdots ,\va_{2N} \}=\{\va_{j}\}$.
Then we have
\begin{align}\label{eq:production of jH and eta GP}
\int_{\T} \la\Hess(\eta)\vj(H),\J\vj(H) \ra d \vx
=&-\nabla \eta(\va_{j})\cdot(\J\nabla_{\va_{j}}\W(\va)).
\end{align}
\end{Lem}
\begin{proof}
Denote
\begin{equation}
  \vj_0=-\J\nabla\Psi.
 \end{equation}
Then \eqref{eq:def of H} and \eqref{eq:def of Phi} imply
\begin{equation}\label{jHwpi}
  \vj(H)=\vj_0+2\pi\vq(\va).
\end{equation}
Plugging \eqref{jHwpi} into \blackmark{ the left hand-side of} \eqref{eq:production of jH and eta GP}, we get
\begin{align}\label{eq:firstestimate of j}
  &\int_{\T} \la\Hess(\eta)\vj(H),\J\vj(H)\ra d \vx\nn
  &\quad=\int_{\T}\la\Hess(\eta)(\vj_0+2\pi\vq(\va)),\J(\vj_0+2\pi\vq(\va))\ra d \vx\nn
  &\quad=\int_{\T}\la\Hess(\eta)\vj_0,\J\vj_0\ra d \vx+2\pi\int_{\T}\la\Hess(\eta)\vj_0,\J\vq(\va)\ra d \vx\nn&\qquad+2\pi\int_{\T}\la\Hess(\eta)\vq(\va),\J\vj_0\ra d \vx+4\pi^2\int_{\T}\la\Hess(\eta)\vq(\va),\J\vq(\va)\ra d \vx.
\end{align}
For the first term of the \blackmark{right hand-side of} \eqref{eq:firstestimate of j}, Lemma 2.3.1 in \cite{RN8} implies that
\begin{equation}\label{eq:production of cunrrent of H}
\int_{\T}\la \Hess(\eta)\vj_0,\J\vj_0\ra d \vx=d_{j} \nabla \eta(\va_{j})\cdot\left(2 \pi\sum_{1\le k\le 2N, \, k\ne j}d_k(\J\nabla F(\va_{j}-\va_k)) \right).
\end{equation}
Applying integration by parts to the second, third and fourth terms of \blackmark{the right hand-side of} \eqref{eq:firstestimate of j}, we have
\begin{align}
\int_{\T}\la\Hess(\eta)\vq(\va),\J\vq(\va) \ra d \vx&=0,\label{eq:twoconst}\\
\int_{\T}\la\Hess(\eta)\vj_0,\J\vq(\va) \ra d \vx&=-\int_{\T}\nabla\cdot \vj_0\la\nabla \eta,\J\vq(\va) \ra d \vx=0,\label{eq:kk}\\
\int_{\T}\la\Hess(\eta)\vq(\va),\J\vj_0 \ra d \vx&=-\int_{\T}\nabla\cdot( \J\vj_0)\la\nabla \eta,\vq(\va) \ra d \vx=-\int_{\T}\vq(\va)\cdot \nabla\eta \nabla\cdot(\J\vj_0)d\vx\nn&=-2\pi\int_{\T}\vq(\va)\cdot \nabla\eta\sum_{k=1}^{2N}d_k\delta_{\va_k}d\vx=-2\pi\vq(\va)\cdot \left(d_{j}\nabla\eta(\va_{j})\right).\label{eq:kl}
\end{align}
In the above equalities,  \blackmark{we have used} $\supp (\eta)\cap \{ \va_1,\cdots ,\va_{2N} \}=\{\va_{j}\}$ and
\begin{equation}
  \nabla\cdot\vj_0=0,\quad \nabla\cdot(\J\vj_0)=2\pi\sum_{j=1}^{2N}d_j\delta_{\va_j},
\end{equation}
which are direct corollaries of \eqref{jHwpi}, \eqref{eq:div of jH} and \eqref{jvevJv}.

Substituting \eqref{eq:production of cunrrent of H}-\eqref{eq:kl} into \eqref{eq:firstestimate of j}, we obtain
\begin{equation}\label{eq:lllajdsf}
  \int_{\T} \la\Hess(\eta)\vj(H),\J\vj(H)\ra d \vx=-2\pi d_{j}\nabla\eta(\va_{j})\cdot\left(-\sum_{1\le k\le 2N, \, k\ne j}d_k(\J\nabla F(\va_{j}-\va_k)) +2\pi\vq(\va)\right).
\end{equation}
Taking the gradient of $W_{\T}$ defined by \blackmark{\eqref{eq:define of W} and \eqref{eq:define of WQ0}} with respect to $\va_j$, we have
\begin{equation}\label{eq:Derivate of W}
  \nabla_{\va_j}\W(\va)=-2 \pi\sum_{1\le k\le 2N,\, k\ne j}d_kd_j\nabla F(\va_j-\va_k)-4 \pi^2 d_j\J\vq(\va).
\end{equation}
Then substituting \eqref{eq:Derivate of W} into \eqref{eq:lllajdsf}, we obtain \eqref{eq:production of jH and eta GP}.
\end{proof}

\section{Vortex paths and current of the NLSE \eqref{eq:maineq}}\label{sec:exist}

In this section, we will first  derive local existence of $2N$ vortex paths by proving
\begin{equation}\label{eq:congvergence of J(u(t))}
    J(u^\varepsilon(t))\to \pi\sum_{j=1}^{2N}d_j \delta_{\vb_j(t)}\ \text{in}\ \M(\T)\quad\text{for some Lipschitz paths}\ \vb_j\text{'s},
\end{equation}
 then prove the convergence of $\vj(u^\varepsilon)$ and $\frac{\vj(u^\varepsilon)}{|u^\varepsilon|}$, and finally give some estimates of $L^2$-norm of  $\nabla |u^\varepsilon|$ and $\frac{\vj(u^\varepsilon)}{|u^\varepsilon|}-\vj(u_*)$ with
  \begin{equation}\label{eq:def of u*}
    u_*(\vx,t)=H(\vx;\vb(t)),
    \end{equation}
    where $\vb:=\vb(t)=(\vb_1(t),\cdots,\vb_{2N}(t))^T$ and $H(\vx;\vb(t))$ is the canonical harmonic map given by \eqref{eq:def of H2}. Similar to the proofs in \cite{RN8,RN2}, the proof of our main result relies on these results.

\subsection{Local existence of vortex paths}

\begin{Lem}[Local existence of vortex paths as $\varepsilon\to 0$]\label{thm:existence and regularity}
If the initial data $u_0^\varepsilon$ satisfies \eqref{con:convergence of J varphi} and \eqref{Con:strong bd of energy}, then the NLSE \eqref{eq:maineq} with \eqref{init}
has a weak solution $u^\varepsilon(\vx,t)\in H^1(\T\times\R)$ for each $\varepsilon>0$. Moreover, \blackmark{after passing to a subsequence $\varepsilon\to 0^+$}, still denoted by $\varepsilon$, there exists a $T>0$ and $2N$ Lipschitz paths $\vb_j:[0,T)\to \T$ with $\vb_j(0)=\va_{j}^0$ for $j=1,2,\cdots,2N$, such that
\eqref{eq:congvergence of J(u(t))} holds for all $t\in [0,T)$. Moreover, $T$ satisfies
\begin{equation}\label{eq:def of T}
    T=\inf \{t>0|\, \exists 1\le j<k\le 2N\ \text{such that}\  |\vb_j(t)-\vb_k(t) |= 0 \}.
\end{equation}

\end{Lem}

\begin{proof}
For the local well-posedness of \blackmark{Problem \eqref{eq:maineq}} with \eqref{init} for each $\varepsilon>0$, one can see \cite{RN4} for details and thus they are omitted here for brevity. The proofs of the existence of $\vb_j$ and \eqref{eq:def of T} are essentially the same as the proof of Theorem 1.4.1 in \cite{RN8}, since the proof of Theorem 1.4.1 in \cite{RN8} does not depend on the vanishing momentum assumption $\bvec{Q}_0=\bvec{0}$.
\end{proof}

\subsection{Convergence of current density}

\begin{Lem}\label{claim:convergence of j/|u|}
  Assume $u^{\varepsilon}, u^\varepsilon_0, \vb_j$ are \blackmark{the same as in} Lemma \ref{thm:existence and regularity}. Then there exists a \blackmark{subsequence} $\varepsilon\to 0$, which is still denoted by $\varepsilon$, such that for any $T_0<T$
    \begin{equation}\label{eq:convergence of j(u(t))}
  \vj(u^\varepsilon)\wto \vj(u_*)\ \text{in}\ L^1(\T\times [0,T_0]),
\end{equation}
and
\begin{equation}\label{eq:convergence of j/||}
\frac{\vj(u^{\varepsilon})}{|u^{\varepsilon}|}\wto \vj(u_*)\ \text{in}\ L^2_{\rm loc}(\blackmark{\T_*(\vb(t))}\times[0,T_0]),
\end{equation}
where $\T_*(\vb(t))$ is defined by replacing $\va$ by $\vb(t)$ in \eqref{eq:def of T_*a}.
\end{Lem}

\begin{proof}
We first prove \eqref{eq:convergence of j(u(t))}.
Lemma \ref{thm:existence and regularity} implies that for $\varepsilon$ small enough, we have for $t\le T_0$,
\begin{equation}
  \left\|\pi\sum_{j=1}^{2N}d_j\delta_{\vb_j(t)}-J(u^\varepsilon(t))\right\|_{W^{-1,1}(\T)}\le \frac{\pi r_{_{\vb}}}{200},
\end{equation}
where
\[
r_{_{\vb}}=\frac{1}{4}\min\bigl\{|\vb_j(t)-\vb_k(t)|1\le j<k\le 2N,0\le t\le T_0 \bigr\}.
\]
The energy conservation \eqref{eq:conserve E}, the energy bound of initial data given in \eqref{Con:strong bd of energy} and the definition of \blackmark{$W^\varepsilon_{\T}$  \eqref{eq:def of We}} imply that there exists a positive constant $C$ such that
\begin{equation}\label{eq:weak bound of E}
  E(u^\varepsilon(t))=E(u^\varepsilon_0)\le 2N\pi\log\frac{1}{\varepsilon}+C.
\end{equation}
Then, it follows from (1.4.26) in \cite{RN8} that  for any $0<\rho<r_{_{\vb}}$,
\begin{equation}\label{eq:bd of int|j|}
\|\vj(u^\varepsilon)\|_{L^1(\T\times [0,T_0])}\le C,\qquad   \|\nabla u^\varepsilon\|_{L^2(\T_\rho(\vb(t)))}\le C,
\end{equation}
which implies that there exists a subsequence of $u^\varepsilon$ and $\vj_*\in L^1(\T\times[0,T_0])$ such that
\begin{equation}\label{eq:con of jue to j*}
\vj(u^\varepsilon)\rightharpoonup \vj_*\  \text{in}\  L^1(\T\times [0,T_0]).
\end{equation}
In addition, the mass conservation \eqref{eq:conserve M} and the energy bound \eqref{eq:weak bound of E} yield
\begin{equation}\label{eq:conservation of L2 norm}
    \frac{\p}{\p t}\frac{|u^\varepsilon|^2}{2}=\Div \vj(u^\varepsilon),\qquad \left\||u^\varepsilon|^2-1\right\|^2_{L^2(\T)}\le C \varepsilon^2|\log \varepsilon|.
\end{equation}
For any $\psi\in C_0^\infty(\T\times[0,T_0])$, noting \eqref{eq:conservation of L2 norm}, we have
\begin{align}\label{eq:int of jue}
\left|\int_0^{T_0}\int_{\T}\vj(u^\varepsilon)\cdot \nabla \psi d \vx\ud t\right|
=&\left|-\int_0^{T_0}\int_{\T}\psi\Div \vj(u^\varepsilon)d \vx\ud t\right|= \left|-\int_0^{T_0}\int_{\T}\psi\p_t\frac{|u^\varepsilon|^2-1}{2}d\vx dt\right|
\nn=&\left|\int_0^{T_0}\int_{\T} \p_t\psi\frac{|u^\varepsilon|^2-1}{2}d \vx\ud t\right|\le C\varepsilon\sqrt{|\log \varepsilon|}\|\p_t\psi\|_{L^2(\T\times[0,T_0])}.
\end{align}
Letting $\varepsilon\to 0$ on both sides of \eqref{eq:int of jue}, we obtain
\[\int_0^{T_0}\int_{\T}\vj_*\cdot\nabla \psi d\vx dt=0,\quad \text{which implies} \ \int_0^{T_0}\int_{\T}\psi\nabla\cdot\vj_* d\vx dt=0.\] Since $\psi$ is arbitrary, we have
\begin{equation}\label{eq:div of u*}
\Div \vj_*=0.
\end{equation}
Similarly, we can prove
\begin{equation}\label{eq:converge of Jue}
\nabla\cdot(\J \vj_*(\vx,t))=\lim_{\varepsilon\to0^+}\nabla\cdot(\J\vj(u^\varepsilon(\vx,t)))=
2\lim_{\varepsilon\to0^+}J(u^\varepsilon(\vx,t))=2 \pi\sum_{k=1}^{2N} d_k \delta _{\vb_k(t)}(\vx).
\end{equation}
Similar to \eqref{eq:relation between Q and q}, \eqref{eq:converge of Jue} implies
\begin{equation}\label{eq:int of ju*}
  \int_{\T}\vj_*(\vx,t)d\vx=\lim_{\varepsilon\to 0}\int_{\T}\vj(u^\varepsilon(\vx,t))d\vx=2\pi\J\sum_{j=1}^{2N}d_j\vb_j(t)=2\pi\vq(\vb(t)).
\end{equation}
Combining \eqref{eq:def of u*}, \eqref{eq:div of jH} and \eqref{eq:int of ju*}, we have
\begin{equation}\label{eq:eq of int j*and int j}
  \int_{\T}\vj(u_*(\vx,t))d\vx=2\pi\vq(\vb(t))=\int_{\T}\vj_*(\vx,t)d\vx.
\end{equation}
Define $\bvec{V}=\vj_*-\vj(u_*)$. Combining \eqref{eq:div of jH}, \eqref{eq:div of u*} and \eqref{eq:converge of Jue}, we see that
\[\Div \bvec{V}=\Div \vj_*-\Div \vj(u_*)=0,\qquad    \nabla\cdot(\J\bvec{V})=\nabla\cdot(\J\vj_*)-\nabla\cdot(\J \vj(u_*))=0.
\]
Thus, $\bvec{V}(\vx,t)=\bvec{g}(t)$ for some $\bvec{g}:[0,T_0]\to \mC$ and $\bvec{g}\in L^1([0,T_0])$ since $\vj_*,\vj(u_*)\in L^1(\T\times[0,T_0])$.   Then, \eqref{eq:eq of int j*and int j} implies that for any $t_1,t_2\in[0,T_0]$,
\begin{align*}
\int_{t_1}^{t_2}\bvec{g}(t)\ud t=\int_{t_1}^{t_2}\int_{\T} \bvec{V}(\vx,t)d \vx\ud t=\int_{t_1}^{t_2}\int_{\T}(\vj_*(\vx,t)-\vj(u_*(\vx,t)))d \vx\ud t=\bvec{0},
\end{align*}
which implies $\bvec{V}(\vx,t)=\bvec{g}(t)=\bvec{0}$. Hence, $\vj_*-\vj(u_*)=\bvec{V}=\bvec{0}$, which together with \eqref{eq:con of jue to j*} implies \eqref{eq:convergence of j(u(t))}.

Then we give the proof of \eqref{eq:convergence of j/||}. Noting
\begin{equation}
  |\nabla u^\varepsilon|^2=\bigl|\nabla |u^\varepsilon|\,\bigr|^2+\left|\frac{\vj(u^\varepsilon)}{|u^\varepsilon|} \right|^2\ge\left|\frac{\vj(u^\varepsilon)}{|u^\varepsilon|} \right|^2,
\end{equation}
\eqref{eq:bd of int|j|} implies  that for any $0<\rho<r_{_{\vb}}$,
\begin{equation}
  \int_0^{T_0}\int_{\T_\rho(\vb(t))}\left|\frac{\vj(u^\varepsilon)}{|u^\varepsilon|}\right|^2d\vx dt\le\int_0^{T_0}\int_{\T_\rho(\vb(t))}\left|\nabla u^\varepsilon\right|^2d\vx dt\le C\,T_0,
\end{equation}
i.e. $\frac{\vj(u^\varepsilon)}{|u^\varepsilon|}$ is uniformly bounded in \blackmark{$L^2(\T_{\rho}(\vb(t))\times[0,T_0])$}. Hence, there exists a function $\,\widetilde{\vj}_*\in L^2(\T_{\rho}(\vb(t))\times[0,T_0])$ such that up to a subsequence
\begin{equation}\label{eq:convergence of vj/|| t_rho}
  \frac{\vj(u^\varepsilon)}{|u^\varepsilon|}\wto \widetilde{\vj}_*, \ \text{in}\ L^2(\blackmark{\T_{\rho}(\vb(t))}\times[0,T_0]).
\end{equation}
Then \eqref{eq:conservation of L2 norm} implies $|u^\varepsilon|\to 1$ in $L^2(\T\times[0,T_0])$. As a result, $\vj(u^{\varepsilon})=|u^{\varepsilon}|\frac{\vj(u^{\varepsilon})}{|u^{\varepsilon} |}$ converges weakly to $\widetilde{\vj}_*$ in $L^1(\T_{\rho}(\vb(t))\times[0,T_0])$. Hence, $\widetilde{\vj}_*=\vj(u_*)$ by \eqref{eq:convergence of j(u(t))}, which together with \eqref{eq:convergence of vj/|| t_rho} implies that $\frac{\vj(u^\varepsilon)}{|u^\varepsilon|}\wto\vj(u_*)$ in $L^2(\T_\rho(\vb(t))\times[0,T_0])$. Noting that $0<\rho<r_{_{\vb}}$ is arbitrary, we obtain \eqref{eq:convergence of j/||}.
\end{proof}

\subsection[Estimates of L2-norms]{Estimates of \texorpdfstring{$L^2$}{}-norm of \texorpdfstring{$\nabla |u^\varepsilon|$}{}  and  \texorpdfstring{$\frac{\vj(u^\varepsilon)}{|u^\varepsilon|}-\vj(u_*)$}{} }

\begin{Lem}\label{lem:estimate of vj-vj}
Assume that $u^\varepsilon,u_0^\varepsilon,\vb_j,u_*$ are the same as in Lemma \ref{claim:convergence of j/|u|}, and $0<\rho<<1$. Then there exists a positive constant $C$ such that for any $[t_1,t_2]\subset[0,T)$, we have
\begin{align}\label{eq:est of e|| j/||-j}
  &\limsup_{\varepsilon\to 0}\int_{t_1}^{t_2}\int_{\T_\rho(\vb(s))}\left[e(|u^\varepsilon(s)|)+
  \left|\frac{\vj(u^\varepsilon(s))}{|u^\varepsilon(s)|}-\vj(u_*(s))
  \right|^2\right]d\vx ds\nonumber\\
  &\ \ \le C\int_{t_1}^{t_2}|\W(\va^0)-\W(\vb(s))|ds.
\end{align}
\end{Lem}
\begin{proof}
By Lemma 3 in \cite{RN23}, we have  for $0<\varepsilon\le \rho$,
\begin{equation}\label{eq:lower bound of e on B}
  \int_{B_\rho(\vb_j(s))}e(u^\varepsilon)d\vx-\left(\gamma+\pi\log\frac{\rho}{\varepsilon}\right)\ge -\Sigma_1(\varepsilon),
\end{equation}
where
\begin{equation}
  \Sigma_1(\varepsilon)=C\frac{\varepsilon}{\rho}\sqrt{\log\frac{\rho}{\varepsilon}}+\frac{C}{\rho}\|J(u)-\pi\delta_{\vb_j(s)}\|_{\M(B_\rho(\vb_j(s)))}.
\end{equation}
Combining \eqref{eq:lower bound of e on B}, \eqref{eq:conserve E} and \eqref{Con:strong bd of energy}, we obtain
\begin{align}\label{eq:bd of e except B}
\int_{\T_\rho(\vb(s))}e(u^\varepsilon)d\vx=&E(u^\varepsilon(s))-
\sum_{j=1}^{2N}\int_{B_\rho(\vb_j(s))}e(u^\varepsilon(s))d\vx\nonumber\\
\le&E(u^\varepsilon_0)-2N\left(\gamma+\pi\log\frac{\rho}{\varepsilon}\right)
+\Sigma_1(\varepsilon)\nonumber\\
  \le&\W^\varepsilon(\va^0)-2N\left(\gamma+\pi\log\frac{\rho}
  {\varepsilon}\right)+\Sigma_1(\varepsilon)\nonumber\\
  =&2N\pi\log\frac{1}{\rho}+\W(\va^0)+\Sigma_1(\varepsilon).
\end{align}
Combining \eqref{eq:est of int of |jH|^2} and \eqref{eq:bd of e except B}, we get
\begin{equation}\label{eq:distance bt ue u*}
  \int_{\T_\rho(\vb(s))}(e(u^\varepsilon)-e(u_*))d\vx\le \W(\va^0)-\W(\vb(s))+O(\rho^2)+\Sigma_1(\varepsilon).
\end{equation}
Noting that
\begin{equation}
  e(u^\varepsilon)-e(u_*)=e(|u^\varepsilon|)+\frac{1}{2}\left|\frac{\vj(u^\varepsilon)}{|u^\varepsilon|}-\vj(u_*) \right|^2+\vj(u_*)\cdot\left(\frac{\vj(u^\varepsilon)}{|u^\varepsilon|}-\vj(u_*)\right),
\end{equation}
\eqref{eq:distance bt ue u*} implies
\begin{align}\label{eq:short time estimate}
  &\int_{\T_\rho(\vb(s))}\left(e(|u^\varepsilon|)+\frac{1}{2}
  \left|\frac{\vj(u^\varepsilon)}{|u^\varepsilon|}-\vj(u_*) \right|^2\right)d\vx\nn
  &\quad\le |\W(\va^0)-\W(\vb(s))|+O(\rho^2)+\Sigma_1(\varepsilon)-\int_{\T_\rho(\vb(s))}\vj(u_*)\cdot\left(\frac{\vj(u^\varepsilon)}{|u^\varepsilon|}-\vj(u_*)\right)d\vx.
\end{align}
Integrating  both sides of \eqref{eq:short time estimate}  with respect to $s$ over $[t_1,t_2]$, we obtain
\begin{align}\label{intttph}
  &\int_{t_1}^{t_2}\int_{\T_\rho(\vb(s))}\left(e(|u^\varepsilon|)+\frac{1}{2}\left|\frac{\vj(u^\varepsilon)}{|u^\varepsilon|}-\vj(u_*) \right|^2\right)d\vx ds \nn&\quad\le \int_{t_1}^{t_2}|\W(\va^0)-\W(\vb(s))|ds+O(\rho^2)+\blackmark{(t_2-t_1)\Sigma_1(\varepsilon)}\nn
  &\qquad-\int_{t_1}^{t_2}\int_{\T_\rho(\vb(s))}\vj(u_*)\cdot\left(\frac{\vj(u^\varepsilon)}
  {|u^\varepsilon|}-\vj(u_*)\right)d\vx ds.
\end{align}
By \eqref{eq:convergence of j/||}, letting $\varepsilon\to 0$ on both sides of \eqref{intttph},  we  obtain
\begin{align}
 & \limsup_{\varepsilon\to 0}\int_{t_1}^{t_2}\int_{\T_\rho(\vb(s))}\left(e(|u^\varepsilon|)
 +\frac{1}{2}\left|\frac{\vj(u^\varepsilon)}{|u^\varepsilon|}-\vj(u_*) \right|^2\right)d\vx ds\nn&\quad
 \le \int_{t_1}^{t_2}|\W(\va^0)-\W(\vb(s))|ds+O(\rho^2).
\end{align}
Here, \blackmark{we have used}
\begin{equation}\label{eq:l2ju}
  \vj(u_*)\in L^2_{loc}(\T_*(\vb(t))\times[t_1,t_2]),
\end{equation}
which is a direct corollary of \eqref{eq:def of u*}, \eqref{eq:def of H} and $F\in C^\infty_{loc}(\T\setminus\{\bvec{0}\})$.

Since the estimate above works for any $\rho'<\rho$ and $\T_{\rho'}(\vb(s))\supset\T_\rho(\vb(s))$, we have
\begin{align}
  &\limsup_{\varepsilon\to 0}\int_{t_1}^{t_2}\int_{\T_\rho(\vb(s))}\left(e(|u^\varepsilon|)
  +\frac{1}{2}\left|\frac{\vj(u^\varepsilon)}{|u^\varepsilon|}-\vj(u_*) \right|^2\right)d\vx ds\nn& \quad\le\limsup_{\varepsilon\to 0}\int_{t_1}^{t_2}\int_{\T_{\rho'}(\vb(s))}\left(e(|u^\varepsilon|)+\frac{1}{2}\left|\frac{\vj(u^\varepsilon)}{|u^\varepsilon|}-\vj(u_*) \right|^2\right)d\vx ds\nn &\quad\le \int_{t_1}^{t_2}|\W(\va^0)-\W(\vb(s))|ds+O((\rho')^2).
\end{align}
Letting $\rho'\to0$, we obtain \eqref{eq:est of e|| j/||-j}.
\end{proof}

\section{Proof of our main result and its extension }\label{sec:proof main}

In this section, we prove the main result Theorem \ref{thm:dynamics} and then discuss
its extension to \blackmark{torus} with arbitrary length and width.

\subsection{Proof of our main result}

\begin{proof}[Proof of Theorem \ref{thm:dynamics}]
Recall $\va$ is the solution of \eqref{eq:dynamics} and $\vb,T$ were obtained in Lemma \ref{thm:existence and regularity}.
Define
\begin{align}\label{eq:define of T'}
 &T_1=\inf \{t>0|\, \exists 1\le j<k\le 2N\ \text{such that}\  |\va_j(t)-\va_k(t) |= 0 \}, \qquad T_2:=\min\{T,T_1\},\\
 \label{eq:def of zeta}
 &\zeta(t):=\sum_{j=1}^{2N}|\vb_j(t)-\va_j(t)|, \qquad t\ge0.
\end{align}

For any $\widetilde{T}<T_2$, Lemma \ref{thm:existence and regularity} and the definition of $\va$ \eqref{eq:dynamics} imply that both $\vb(t)$ and $\va(t)$ are Lipschitz on $[0,\widetilde{T}]$. Hence, there exists a constant $C_*$ such that
\begin{equation}
  |\va(t)-\va(0)|+|\vb(t)-\vb(0)|\le C_*t,
\end{equation}
which together with \eqref{eq:def of zeta} implies that
\begin{equation}\label{eq:lip of zeta}
  |\zeta(t)-\zeta(0)|\le 2NC_*t.
\end{equation}
Lemma \ref{thm:existence and regularity} and the definition of $\va$ \eqref{eq:dynamics} imply that $\vb(0)=\va^0,\va(0)=\va^0$. Hence, \eqref{eq:def of zeta} implies
\begin{equation}\label{eq:zeta(0)=0}
  \zeta(0)=\sum_{j=1}^{2N}|\vb_j(0)-\va_j(0)|=\sum_{j=1}^{2N}|\va_j^0-\va_j^0|=0.
\end{equation}
Then, combining \eqref{eq:lip of zeta} and \eqref{eq:zeta(0)=0} we have
\begin{equation}\label{zetatr4}
\zeta(t)\le \frac{r}{4}, \qquad \text{for}\  0\le t\le \tau_0:=\frac{r}{8NC_*},
\end{equation}
with
\begin{equation}\label{requalab}
r=\frac{1}{4}\min\left\{|\va_j(t)-\va_k(t) |,|\vb_j(t)-\vb_k(t) |\,\left|1\le j<k\le 2N,0\le t\le \widetilde{T} \right.\right\}.
\end{equation}

We first  prove $\vb(t)=\va(t)$, which is equivalent to $\zeta(t)=0$ on $[0,\tau_0]$.

Substituting \eqref{eq:dynamics} into the derivative of $\zeta(t)$, we have
\begin{align}\label{eq:splitting of dz}
    \dot{\zeta}(t)&\le \sum_{j=1}^{2N}\left|\dot{\va}_j(t)-\dot{\vb}_j(t) \right|=\sum_{j=1}^{2N}\left|-d_j\frac{1}{\pi}\J\nabla_{\va_j} \W(\va(t)) -\dot{\vb}_j(t) \right|\nn
    &\le \sum_{j=1}^{2N}\left|-d_j\frac{1}{\pi}\J\nabla_{\va_j} \W(\va(t))+d_j\frac{1}{\pi}\J\nabla_{\vb_j} \W(\vb(t))\right|+\sum_{j=1}^{2N}\left|-d_j\frac{1}{\pi}\J\nabla_{\vb_j} \W(\vb(t)) -\dot{\vb}_j(t) \right|\nn
    &= \sum_{j=1}^{2N}B_{j}(t)+\sum_{j=1}^{2N}A_j(t),
\end{align}
where
\begin{align}
A_j(t)&=\left|\dot{\vb}_j(t)+d_j\frac{1}{\pi}\J\nabla_{\vb_j} \W(\vb(t)) \right|,\\
B_{j}(t)&=\left|-d_j\frac{1}{\pi}\J\nabla_{\va_j} \W(\va(t))+d_j\frac{1}{\pi}\J\nabla_{\vb_j} \W(\vb(t)) \right|.\label{eq:define of B}
\end{align}

Substituting \eqref{eq:Derivate of W} into the definition of $B_j$ in \eqref{eq:define of B}, we have
\begin{align}\label{eq:first etimate of B}
B_j(t)\le&2\left|\sum_{1\le k\le 2N,\, k\ne j}d_k\J\nabla F(\va_j(t)-\va_k(t))-\sum_{1\le k\le 2N,\, k\ne j}d_k\J\nabla F(\vb_j(t)-\vb_k(t))\right|\nn
&+\left|-4 \pi\J\sum_{k=1}^{2N}d_k\va_k(t)+4 \pi \J\sum_{k=1}^{2N}d_k\vb_k(t) \right|\nn
\le&2\sum_{\begin{subarray}{c}1\le k\le 2N,k\ne j\end{subarray}}|d_k\J\nabla F(\va_j(t)-\va_k(t))-d_k\J\nabla F(\vb_j(t)-\vb_k(t))|+8\pi N\zeta(t)\nn
\le&2\zeta(t)\sum_{\begin{subarray}{c}1\le k\le 2N,k\ne j\end{subarray}}\|F\|_{C^2(B_r(\vb_j(t)-\vb_k(t)))}+8\pi N\zeta(t).
\end{align}
In the above, \blackmark{we have used} the following inequality by noting \eqref{eq:def of zeta} and
\eqref{zetatr4}
\[|(\va_j-\va_k)-(\vb_j-\vb_k)|\le |\vb_j-\va_j|+|\vb_k-\va_k|\le \zeta<r,
\qquad 1\le j\ne k\le 2N.
\]
Since $B_r(\vb_j(t)-\vb_k(t))\subset\T\setminus B_{3r}(\bvec{0})$, there exists a positive constant $C$ such that $\|F\|_{C^2(B_r(\vb_j(t)-\vb_k(t)))}\le C$. Hence, \eqref{eq:first etimate of B} implies
\begin{equation}\label{eq:estimate of B}
B_{j}(t)\le (4NC+8\pi N)\zeta(t).
\end{equation}

For $A_j(t)$, we can find a smooth function $\eta\in C_0(B_r(\vb_j(t)))$ satisfying
\[\eta=\eta(\vx)= \bvec{\nu}\cdot \vx\ \text{in}\ B_{3r/4}(\vb_i(t)),
\]
where $\bvec{\nu}\in {\mathbb S}^1$ satisfies
\begin{equation}\label{eq:Ajand vector}
A_j(t)=d_j\, \bvec{\nu}\cdot\left(\dot{\vb_j}(t)+d_j\frac{1}{\pi}\J\nabla_{\vb_j} \W(\vb(t)) \right).
\end{equation}
By  (2.1.9) in \cite{RN8} and \eqref{eq:congvergence of J(u(t))}, we know
\begin{align}\label{eq:mainequality}
d_j \, \bvec{\nu}\cdot \dot{\vb}_j(t)&=\lim_{h\to 0}d_j \,\bvec{\nu}\cdot \frac{1}{h}(\vb_j(t+h)-\vb_j(t))\nn
&=\lim_{h\to 0}\lim_{\varepsilon\to 0}\frac{1}{\pi h}\int_{\T}\left[\eta(\vx) J(u^{\varepsilon}(\vx,t+h))-\eta(\vx) J(u^{\varepsilon}(\vx,t))\right]d\vx\nn
 &=\lim_{h\to 0}\lim_{\varepsilon\to 0}\frac{1}{\pi h}\int_{t}^{t+h}\int_{\T}\la\Hess(\eta)\nabla u^\varepsilon,\J\nabla u^\varepsilon \ra d \vx\ud s.
\end{align}
Lemma \ref{lem:production of jH and eta} together with \eqref{eq:def of u*} implies
\begin{align}\label{eq:vector nabla w}
d_j\,\bvec{\nu}\cdot\left(d_j\frac{1}{\pi}\J\nabla_{\vb_j} \W(\vb(t))\right)=&-\frac{1}{\pi}\int_{\T}\la \Hess(\eta)\vj(u_*(t)),\J\vj(u_*(t)) \ra d \vx\nn
=&-\lim_{h\to 0}\frac{1}{\pi h}\int_t^{t+h}\int_{\T}\la \Hess(\eta)\vj(u_*(s)),\J\vj(u_*(s)) \ra d \vx\ud s.
\end{align}
Combining \eqref{eq:Ajand vector}, \eqref{eq:mainequality} and \eqref{eq:vector nabla w}, and noting
\begin{equation}\label{eq:separation of nabla square}
  \begin{cases}
  \la\p_xu^\varepsilon,\p_xu^\varepsilon\ra&=\frac{1}{|u^\varepsilon|^2}(j_1(u^\varepsilon))^2+(\p_x|u^\varepsilon|)^2,\\
  \la\p_xu^\varepsilon,\p_yu^\varepsilon\ra&=\frac{1}{|u^\varepsilon|^2}j_1(u^\varepsilon)j_2(u^\varepsilon)+\p_x|u^\varepsilon|\,\p_y|u^\varepsilon|,\\
  \la\p_yu^\varepsilon,\p_yu^\varepsilon\ra&=\frac{1}{|u^\varepsilon|^2}(j_2(u^\varepsilon))^2+(\p_y|u^\varepsilon|)^2, 
  \end{cases}
\end{equation}
which is a corollary of 
\begin{equation}
  \nabla u^\varepsilon=\frac{u^\varepsilon}{|u^\varepsilon|}\nabla|u^\varepsilon|+\frac{\vj(u^\varepsilon)}{|u^\varepsilon|}\frac{iu^\varepsilon}{|u^\varepsilon|},
\end{equation}
one gets
\begin{align}\label{eq:splitting of A}
    A_j(t)=&\lim_{h\to 0}\lim_{\varepsilon\to 0}\frac{1}{\pi h}\int_t^{t+h}\int_{\T}\bigl[\la\Hess(\eta)\nabla u^{\varepsilon},\J\nabla u^{\varepsilon} \ra-\la \Hess(\eta)\vj(u_*),\J\vj(u_*) \ra\bigr]d \vx\ud s\nonumber\\
    =&L_j(t)+K_j(t)=L_j(t)+K_{j1}(t)+K_{j2}(t)+K_{j3}(t),
\end{align}
where
\begin{align}
    L_j(t)=&\lim_{h\to 0}\lim_{\varepsilon\to 0}\frac{1}{\pi h}\int_t^{t+h}\int_{\T}\la\Hess(\eta)\nabla|u^{\varepsilon}|,\J\nabla|u^{\varepsilon}| \ra d \vx\ud s,\label{eq:define of L}\\
K_j(t)=&\lim_{h\to 0}\lim_{\varepsilon\to 0}\frac{1}{\pi h}\int_t^{t+h}\int_{\T}\left[\la\Hess(\eta)\frac{\vj(u^{\varepsilon})}{|u^{\varepsilon}|},\frac{\J\vj(u^{\varepsilon})}{|u^{\varepsilon}|} \ra-\la \Hess(\eta)\vj(u_*),\J\vj(u_*) \ra\right]d \vx\ud s,\\
K_{j1}(t)=&\lim_{h\to 0}\lim_{\varepsilon\to 0}\frac{1}{\pi h}\int_t^{t+h}\int_{\T}\la\Hess(\eta)\left(\frac{\vj(u^{\varepsilon})}{|u^{\varepsilon}|}-\vj(u_*)\right),\J\left(\frac{\vj(u^{\varepsilon})}{|u^{\varepsilon}|}-\vj(u_*)\right)\ra d \vx\ud s,\label{eq: 1st term of L} \\
K_{j2}(t)=&\lim_{h\to 0}\lim_{\varepsilon\to 0}\frac{1}{\pi h}\int_t^{t+h}\int_{\T}\la\Hess(\eta)\vj(u_*),\J\left(\frac{\vj(u^{\varepsilon})}{|u^{\varepsilon}|}-\vj(u_*)\right)\ra d \vx\ud s,\label{eq: 2nd term of L} \\
K_{j3}(t)=&\lim_{h\to 0}\lim_{\varepsilon\to 0}\frac{1}{\pi h}\int_t^{t+h}\int_{\T}\la\Hess(\eta)\left(\frac{\vj(u^{\varepsilon})}
{|u^{\varepsilon}|}-\vj(u_*)\right),\J\vj(u_*)\ra d \vx\ud s.\label{eq: 3rd term of L}
{}\end{align}
Noting \eqref{eq:l2ju} and $\Hess(\eta)\in C_0(B_r(\vb_j(t))\setminus B_{3r/4}(\vb_j(t)))$, substituting \eqref{eq:convergence of j/||} to \eqref{eq: 2nd term of L} and \eqref{eq: 3rd term of L}, we obtain
\begin{equation}\label{eq:estimate of Kj2 Kj3}
    K_{j2}(t)\equiv 0, \qquad K_{j3}(t)\equiv 0.
\end{equation}
Thus, it only remains to estimate $L_j(t)$ and $K_{j1}(t)$.
By definition \eqref{eq:define of L} and \eqref{eq: 1st term of L}, 
\begin{align*}
|K_{j1}(t)|+ |L_j(t)|\le C\lim_{h\to 0}\lim_{\varepsilon\to 0}\frac{1}{\pi h}\int_t^{t+h}\left(\left\|\frac{\vj(u^\varepsilon)}{|u^\varepsilon|}-\vj(u_*)\right\|_{L^2(\T_\rho(\vb(s)))}+\left\|\nabla |u^\varepsilon|\right\|_{L^2(\T_\rho(\vb(s)))}\right)ds.
\end{align*}
Then by Lemma \ref{lem:estimate of vj-vj}, we have
\begin{align}\label{eq:estimate of Kj1 L}
  |K_{j1}(t)|+ |L_j(t)|\le& C\lim_{h\to 0}\frac{1}{\pi h}\int_t^{t+h}|\W(\va^0)-\W(\vb(s))|ds\le \frac{C}{\pi}|\W(\va^0)-\W(\vb(t))|.
\end{align}
Noting \eqref{eq:define of WQ0} and \eqref{eq:dynamics}, differentiate $\W(\va(t))$ with respect to $t$:
\begin{align}\label{eq:conservation law of W}
\frac{\ud }{\ud t}\W(\va(t))&=\sum_{j=1}^{2N}\nabla_{\va_j}\W(\va(t))\cdot \dot{\va}_j(t)=\nabla_{\va_j}\W(\va)\cdot\left(-\frac{1}{\pi}d_j\J\nabla_{\va_i}\W(\va)\right)=0.
\end{align}
This implies that $\W(\va(t))\equiv \W(\val)$.
As a result, \eqref{eq:estimate of Kj1 L} gives
\begin{equation}\label{eq:final estimate of Kj1 L}
  |K_{j1}(t)|+|L_j(t)|\le \frac{C}{\pi}|\W(\va(t))-\W(\vb(t))|\le C\zeta(t).
\end{equation}

Combining \eqref{eq:splitting of dz}, \eqref{eq:estimate of B}, \eqref{eq:splitting of A}, \eqref{eq:estimate of Kj2 Kj3} and \eqref{eq:final estimate of Kj1 L}, we obtain
\[\dot{\zeta}(t)\le \sum_{j=1}^{2N}B_{j}(t)+\sum_{j=1}^{2N}(L_j(t)+K_{j1}(t)+K_{j2}(t)+K_{j3}(t))\le C \zeta(t),
\]
which implies $\zeta(t)\equiv 0$ on $[0,\tau_0]$ together with  \eqref{eq:zeta(0)=0}. Hence, $\vb(t)= \va(t)$ on $[0,\tau_0]$. Recall the definition of $\tau_0$ in \eqref{requalab}. $\tau_0$ is a constant whenever $\widetilde{T}$ is chosen. Hence we can repeat the above proof on $[\tau_0,2\tau_0],[2\tau_0,3\tau_0]$ and so on. Then we have $\vb(t)=\va(t)$ for any $t\in[0,\widetilde{T}]$. Since $\widetilde{T}<T_2$ is arbitrary, we have $\vb(t)=\va(t)$ for any $t\in[0,T_2)$. Then  \eqref{eq:def of T} and \eqref{eq:define of T'} imply that $T=T'=\widetilde{T}$. Noting \eqref{eq:congvergence of J(u(t))} and the equivalence between \eqref{newODE} and \eqref{eq:dynamics}, $\vb(t)=\va(t)$ \blackmark{both satisfy}  \eqref{eq:con of J} and \eqref{newODE} on $[0,T)$.
\end{proof}

\subsection{Extension to \blackmark{torus} with arbitrary length and width}
We can extend our result to the case of the nonlinear \SE{}  on \blackmark{torus} with arbitrary length and width \blackmark{$\T_{lw}=(\R/l\Z)\times(\R/w\Z)$ with $l>0$ and $w>0$}:
\begin{equation}\label{eq:maineqgt}
    \I \p_tu^\varepsilon(\vx,t)- \Delta u^\varepsilon(\vx,t)+\frac{1}{\varepsilon^2}\left(|u^\varepsilon(\vx,t)|^2-1 \right)u^\varepsilon(\vx,t)=0,\quad  \blackmark{\vx\in\T_{lw}}, \ t>0,
\end{equation}
with initial data
\begin{equation}\label{initgt}
u^\varepsilon(\vx,0)= u_0^\varepsilon(\vx),\quad \blackmark{\vx\in\T_{lw}}.
\end{equation}

Define the renormalized energy \blackmark{$W_{lw}$ on $\T_{lw}$} for $\va=(\va_1,\cdots,\va_{2N})^T\in(\blackmark{\T_{lw}})^{2N}$ with $\va_j\ne \va_k(j\ne k)$ by
\begin{equation}
  \blackmark{W_{lw}}(\va)=-\pi\sum_{1\le k\ne m\le 2N}d_kd_m\blackmark{F_{lw}}(\va_k-\va_m)+\frac{2\pi^2}{\blackmark{lw}}\left|\sum_{k=1}^{2N}d_k\va_k\right|^2,
\end{equation}
with $\blackmark{F_{lw}}$ the solution of
\[
    \Delta \blackmark{F_{lw}}(\vx)=2\pi\left(\delta(\vx)-\frac{1}{\blackmark{lw}} \right)\ \text{on}\ \blackmark{\T_{lw}},\quad \text{with}\ \int_{\blackmark{\T_{lw}}}\blackmark{F_{lw}}(\vx)d\vx=0.
\]
Then we can repeat \blackmark{the proof of Theorem \ref{thm:dynamics}} with some adjustments to prove the following result:

\begin{Coro}[\blackmark{Reduced dynamical law} for the NLSE on \blackmark{torus} with arbitrary length and width]\label{thm:dynamicsgt}
Assume there exist $2N$ distinct points $\va_1^0,\cdots,\va_{2N}^0\in\blackmark{\T_{lw}}$ and $\va^0=(\va_1^0,\cdots,\va_{2N}^0)^T$ such that the initial data $u_0^\varepsilon$ in \eqref{initgt} satisfies
\begin{equation}
  J(u^\varepsilon_0)\stackrel{\varepsilon\to 0^+}{\longrightarrow}\pi\sum_{j=1}^{2N}d_j\delta_{\va_j^0}\ \text{in}\ W^{-1,1}(\blackmark{\T_{lw}}),\qquad\lim_{\varepsilon\to0} \int_{\blackmark{\T_{lw}}}\vj(u^\varepsilon_0(\vx))d\vx=\widehat{\bvec{Q}}_0=2\pi\J\sum_{j=1}^{2N}d_j\va_j^0,
\end{equation}
\begin{equation}\label{Con:strong bd of energygt}
\limsup_{\varepsilon\to 0}\left[\int_{\blackmark{\T_{lw}}}e(u^\varepsilon_0(\vx))d\vx-
2N\left(\pi\log\frac{1}{\varepsilon}+\gamma\right)-\blackmark{W_{lw}}(\val)\right]\le 0.
 \end{equation}
Then there exists a time $\widehat{T}>0$ and $2N$ Lipschitz paths $\va_j:[0,\widehat{T})\to \blackmark{\T_{lw}}$ for $j=1,\cdots,2N$, such that the solution $u^\varepsilon$ of the NLSE \eqref{eq:maineqgt} with
\eqref{initgt} satisfies
\begin{equation}
J(u^\varepsilon(t))\stackrel{\varepsilon\to 0^+}{\longrightarrow} \pi\sum_{j=1}^{2N} d_j \delta_{\va_{j}(t)} \ \text{in}\ \M(\blackmark{\T_{lw}}),
\end{equation}
and $\va_j$ ($1\le j\le 2N$) satisfy \blackmark{the following reduced dynamical law}:
\begin{equation}\label{eq:dynamicsgt}
\dot{\va}_j=-d_j\frac{1}{\pi}\J\nabla_{\va_j}\blackmark{W_{lw}}(\va)=
\J\left(2\sum_{1\le k\le 2N,\, k\ne j}d_k \nabla \blackmark{F_{lw}}(\va_j-\va_k)\right)-\frac{2}{\blackmark{lw}}\bvec{Q}_0 ,\quad  t>0,
\end{equation}
with the initial data $\va_j(0)=\va_j^0$ for $1\le j\le 2N$.
\end{Coro}

\section{Some properties of \blackmark{the reduced dynamical law}}\label{sec:pro of rdl}

In this section, we show some first integrals of \blackmark{the reduced dynamical law} \eqref{newODE} (or \eqref{eq:dynamics}) and present analytical solutions for
several initial setups with symmetry.

\subsection{First Integrals}\label{sec:fi}

Define
\begin{equation}\label{eq:1st integrals}
    \xi(\va):=\frac{1}{4} \sum_{1\le j\ne k\le 2N}d_jd_k|\va_j-\va_k |^2, \qquad \va \in (\T)^{2N}_*.
\end{equation}
Then we have

\begin{Lem}
Let $\va:=\va(t)=(\va_1(t),\ldots,\va_{2N}(t))^T\in (\T)^{2N}_*$ be
the solution of \blackmark{the reduced dynamical law} \eqref{newODE} (or \eqref{eq:dynamics}) with \eqref{initODE}. Then $\vq(\va)$,
$\W(\va)$ and $\xi(\va)$ defined in \eqref{vqva}, \eqref{eq:define of WQ0} and \eqref{eq:1st integrals}, respectively, are three first integrals, i.e.
\begin{eqnarray}\label{vqva2}
&&\vq(\va):=\vq(\va(t))\equiv \vq(\va(0))=\vq(\va^0)=\frac{1}{2\pi}\bvec{Q}_0,\\
\label{wava2}
&&\W(\va):=\W(\va(t))\equiv \W(\va(0))=\W(\va^0),\qquad t\ge0,\\
\label{xiava2}
&&\xi(\va):=\xi(\va(t))\equiv \xi(\va(0))=\xi(\va^0).
\end{eqnarray}
\end{Lem}

\begin{proof}
Differentiating \eqref{vqva} with respect to $t$, noting \eqref{eq:dynamics},
\eqref{d1d2d3}, \eqref{eq:Derivate of W}, and \blackmark{that} $F$ is an even function,  we have
\begin{align*}
\frac{d}{dt}\vq(\va(t))&=2\J\sum_{j=1}^{2N}d_j\J\left(\sum_{1\le k\le 2N,\, k\ne j}d_k\nabla F(\va_j(t)-\va_k(t))-2\pi\sum_{k=1}^{2N}d_k\va_k(t) \right)\\
&=-2\sum_{1\le j\ne k\le 2N}d_jd_k\nabla F(\va_j(t)-\va_k(t))+4\pi\left(\sum_{j=1}^{2N}d_j\right)\left(\sum_{k=1}^{2N} d_k\va_k(t)\right)\\
&=\bvec{0},\qquad t\ge0,
\end{align*}
which immediately implies \eqref{vqva2}. Plugging \eqref{vqva2} into \eqref{eq:dynamics}, we obtain \eqref{newODE} immediately.
From \eqref{eq:conservation law of W},
we get \eqref{wava2}. Finally, combining $\xi(\va)=-\frac{1}{4}|\vq(\va)|^2$
with \eqref{vqva2}, we get \eqref{xiava2}.
\end{proof}

\subsection{Analytical solutions for several initial setups with symmetry}

\begin{Lem}\label{ASODEN1}
When $N=1$ in \eqref{newODE}, the analytical solution of \eqref{newODE} with \eqref{initODE} is given as
\begin{equation}\label{eq:sol of ODE 2v}
\va_1(t)=\va_1^0+\bvec{p}t, \qquad \va_2(t)=\va_2^0+\bvec{p}t,
\end{equation}
with 
\begin{equation}
\blackmark{\bvec{p}=-2\J\nabla F(\va_1^0-\va_2^0)-2\bvec{Q}_0}
\end{equation}
\end{Lem}

\begin{proof}
Clearly, \blackmark{\eqref{eq:sol of ODE 2v} implies} that $\va(t)=(\va_1(t),\va_2(t))^T$ satisfies the initial data of \eqref{newODE} and 
\begin{equation}\label{eq:div of a1-a2}
  \blackmark{\va_1(t)-\va_2(t)\equiv \va_1^0-\va_2^0}, \qquad t\ge0.
\end{equation}
Noting that $F$ is an even function and substituting \eqref{eq:div of a1-a2} into \eqref{newODE}, we have that \blackmark{the right hand-side of} \eqref{newODE} \blackmark{is equal to} ($j=1,2$ respectively)
\begin{equation}
\begin{aligned}
  &2d_2\J \nabla F(\va_1(t)-\va_2(t))-2\bvec{Q}_0=-2\J\nabla F(\va_1^0-\va_2^0)-2\bvec{Q}_0=\bvec{p},\\
  &2 d_1\J\nabla F(\va_2(t)-\va_1(t))-2\bvec{Q}_0=2\J\nabla F(-(\va_1^0-\va_2^0))-2\bvec{Q}_0=-2\J\nabla F(\va_1^0-\va_2^0)-2\bvec{Q}_0=\bvec{p}.
  \end{aligned}
\end{equation}
 Differentiating $\va_j$'s defined by \eqref{eq:sol of ODE 2v} with respect to $t$, we have
\begin{equation}
  \dot{\va}_1(t)=\bvec{p},\qquad \dot{\va}_2(t)=\bvec{p},
\end{equation}
which are equal to the \blackmark{right hand-side of} \eqref{newODE}. By the uniqueness of the solution of \eqref{newODE}, \eqref{eq:sol of ODE 2v} is the solution of \eqref{newODE} when $N=1$.
\end{proof}

\begin{figure}[htp!]
\begin{center}
\begin{tabular}{cccc}
\includegraphics[height=4.5cm]{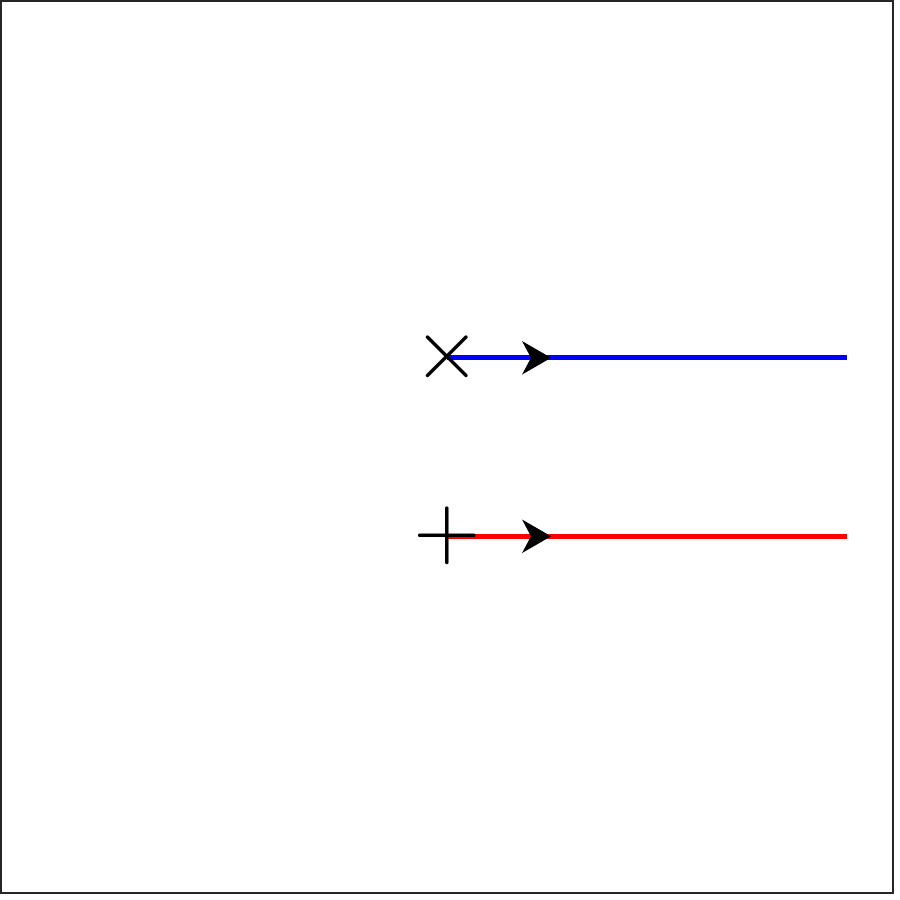}
&\includegraphics[height=4.5cm]{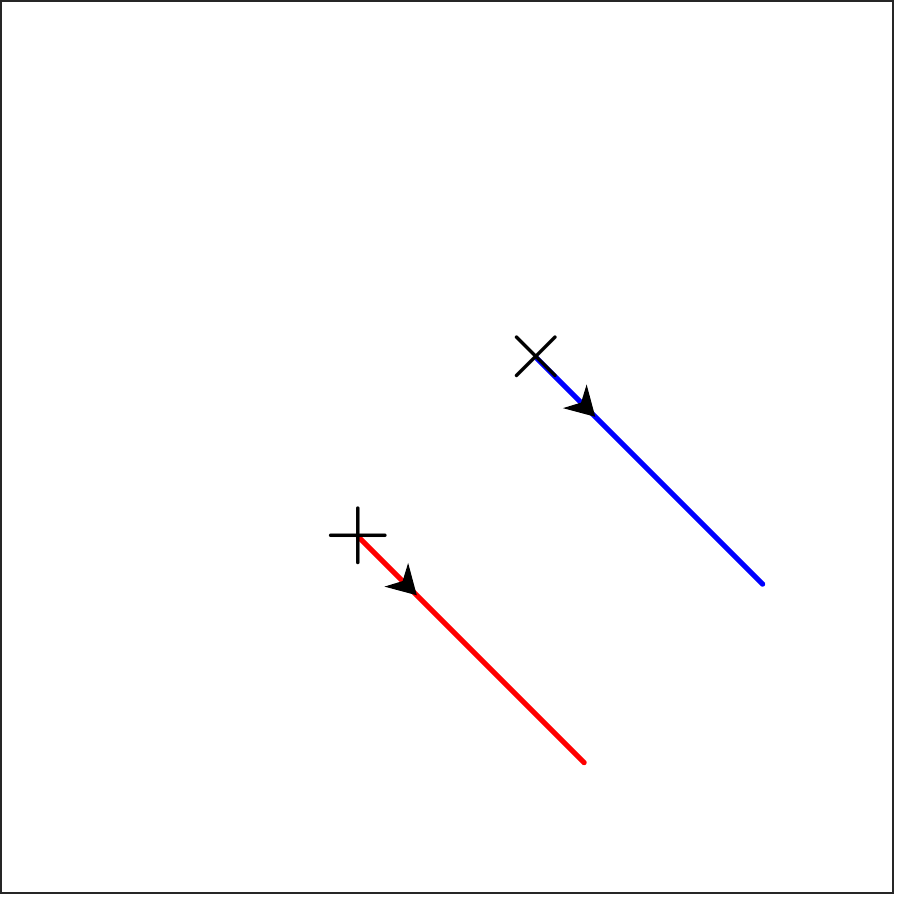}
&\includegraphics[height=4.5cm]{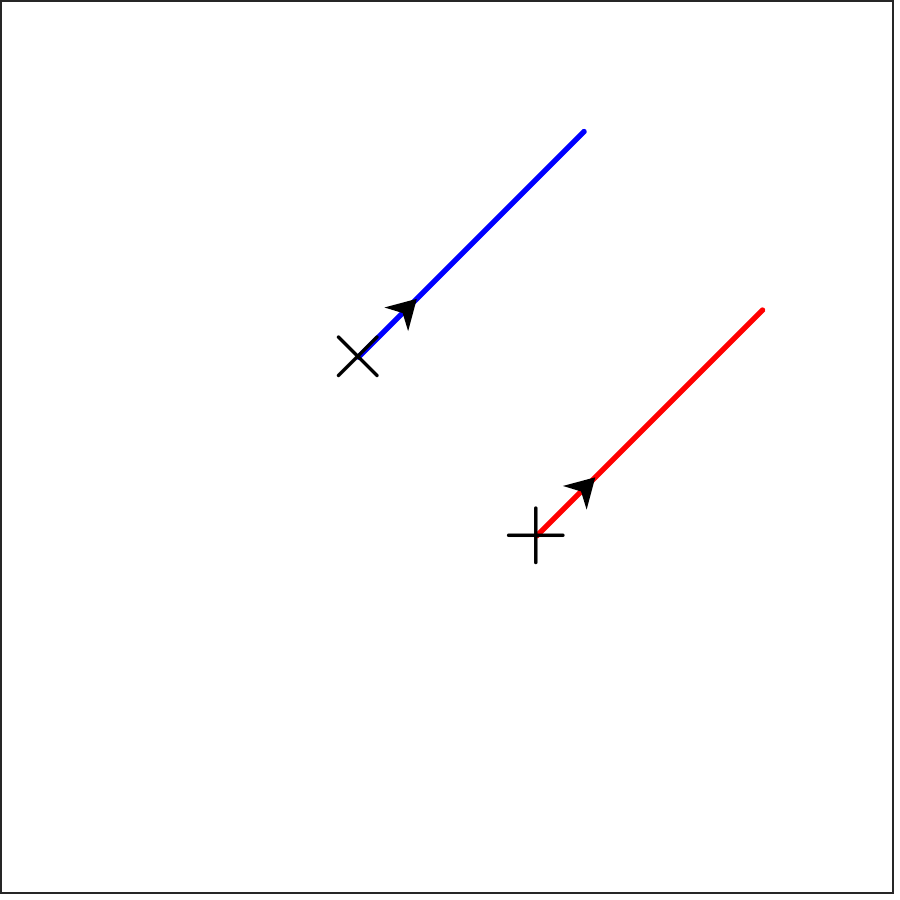}
\end{tabular}
\end{center}
\caption{Trajectories of the reduced dynamics law \eqref{newODE} with $N=1$ and 
different initial datum in \eqref{initODE}, i.e. different $\bvec{Q}_0$, for: (i) 
$\va^0_1=(0.5,0.4)^T$ and $\va^0_2=(0.5,0.6)^T$ with $\bvec{Q}_0=(-0.4\pi,0)^T$ (left), (ii) 
$\va^0_1=(0.4,0.4)^T$ and $\va^0_2=(0.6,0.6)^T$ with $\bvec{Q}_0=(-0.4\pi,0.4\pi)^T$ (middle), and (iii) 
$\va^0_1=(0.6,0.4)^T$ and $\va^0_2=(0.4,0.6)^T$ with $\bvec{Q}_0=(-0.4\pi,-0.4\pi)^T$ (right). Here and in the below, we use $+$ and $\times$  in the pictures to denote vortices with winding numbers $+1$ and $-1$, respectively.
}\label{fig:2vpaths}
\end{figure}

\begin{Lem}\label{lem:2circles}
Take $N=2$ in \eqref{newODE} and the initial data \eqref{initODE} as
\begin{equation}\label{eq:initial of so}
\begin{aligned}
\va^0_1=(0.5,0.5)^T+(\alpha_0,\beta_0)^T,\quad\va^0_2=(0.5,0.5)^T-(\alpha_0,\beta_0)^T,\\ \va_3^0=(0.5,0.5)^T+(\beta_0,\alpha_0)^T,\quad\va^0_4=(0.5,0.5)^T-(\beta_0,\alpha_0)^T,
\end{aligned}
\end{equation}
with $0<\beta_0,\alpha_0<1$ such that $\bvec{Q}_0=\bvec{0}$, then the  analytical solution of \eqref{newODE} with \eqref{initODE} is given as
\begin{equation}\label{eq:symmetry of sol}
\begin{aligned}
\va_1(t)=(0.5,0.5)^T+(\alpha(t),\beta(t))^T, \quad \va_2(t)=(0.5,0.5)^T-(\alpha(t),\beta(t))^T,\\ \va_3(t)=(0.5,0.5)^T+(\beta(t),\alpha(t))^T, \quad \va_4(t)=(0.5,0.5)^T-(\beta(t),\alpha(t))^T,
\end{aligned}
\end{equation}
where ($\alpha(t)$, $\beta(t)$) is the solution of 
\begin{equation}\label{eq:dynamics of alphabeta}
\left\{
\begin{aligned}
  &\dot{\alpha}=2\left(-\p_yF(\alpha-\beta,\beta-\alpha)+\p_yF(2\alpha,2\beta)-\p_yF(\alpha+\beta,\alpha+\beta) \right),\\
  &\dot{\beta}=2(\p_xF(\alpha-\beta,\beta-\alpha)-\p_xF(2\alpha,2\beta)+\p_xF(\alpha+\beta,\alpha+\beta)),
  \end{aligned}
  \right.
\end{equation}
with the initial data
\begin{equation}\label{eq:reduced initial}
  \alpha(0)=\alpha_0,\qquad \beta(0)=\beta_0.
\end{equation}

\end{Lem}

\begin{proof}
By the symmetry of \eqref{eq:define of F}, we \blackmark{have that $F$} satisfies
\begin{equation}
  F(x,y)=F(-x,y)=F(x,-y)=F(y,x).
\end{equation}
\blackmark{Then, owing to the symmetry of} the initial data \eqref{eq:initial of so} and the symmetry of the equation \eqref{newODE},
we can take the ansatz that the solution $\va=(\va_1,\va_2,\va_3,\va_4)^T$ satisfies \eqref{eq:symmetry of sol}. Substituting \eqref{eq:symmetry of sol} into \eqref{newODE} and noting that
\begin{equation}
  \bvec{Q}_0=2\pi\J(\va_1^0+\va_2^0-\va_3^0-\va_4^0)=\bvec{0},
\end{equation}
we have 
\begin{equation}
  \dot{\va}_1=2\J(-\nabla F(\alpha-\beta,\beta-\alpha)+\nabla F(2\alpha,\beta)-\nabla F(\alpha+\beta,\alpha+\beta)).
\end{equation}
Noting $\dot{\va}_1=(\dot{\alpha},\dot{\beta})$, we obtain \eqref{eq:dynamics of alphabeta}.
\end{proof}

\Bluemark{Similar to  Lemma \ref{lem:2circles}, we can prove the following lemmas:
\begin{Lem}\label{SODEN21}
Take $N=2$ in \eqref{newODE} and  the initial data \eqref{initODE} as
\begin{equation}\label{initN21}
\begin{aligned}
&\va^0_1=(0.5,0.5)^T+(\alpha_0,\beta_0)^T,\quad\va^0_2=(0.5,0.5)^T-(\alpha_0,\beta_0)^T,\\ &\va_3^0=(0.5,0.5)^T+(\alpha_0,-\beta_0)^T,\quad\va^0_4=(0.5,0.5)^T-(\alpha_0,-\beta_0)^T,
\end{aligned}
\end{equation}
with $0<\beta_0,\alpha_0<1$ such that $\bvec{Q}_0=\bvec{0}$,
then the  analytical solution of \eqref{newODE} with \eqref{initODE} is given as
\begin{equation}\label{eq:symmetry of sol4cc}
\begin{aligned}
&\va_1(t)=(0.5,0.5)^T+(\alpha(t),\beta(t))^T, \quad \va_2(t)=(0.5,0.5)^T-(\alpha(t),\beta(t))^T,\\ &\va_3(t)=(0.5,0.5)^T+(\alpha(t),-\beta(t))^T, \quad \va_4(t)=(0.5,0.5)^T-(\alpha(t),-\beta(t))^T,
\end{aligned}
\end{equation}
where ($\alpha(t)$,$\beta(t)$) is the solution of 
\begin{equation}\label{eq:dynamics of alphabeta4cc}
\left\{
\begin{aligned}
  &\dot{\alpha}=2\left(\p_yF(2\alpha,2\beta)-\p_yF(0,2\beta) \right),\\
  &\dot{\beta}=2(-\p_xF(2\alpha,2\beta)+\p_xF(2\alpha,0)),
  \end{aligned}
  \right.
\end{equation}
with the initial data \eqref{eq:reduced initial}.
\end{Lem}

\begin{Lem}\label{SODEN22}
Take $N=2$ in \eqref{newODE} and the initial data \eqref{initODE} as
\begin{equation}\label{eq:initial of so4c}
\begin{aligned}
&\va^0_1=(x_0,0.25)^T+(\alpha_0,\beta_0)^T,\quad\va^0_2=(x_0,0.25)^T-(\alpha_0,\beta_0)^T,\\ &\va_3^0=(x_0,0.75)^T+(\alpha_0,-\beta_0)^T,\quad\va^0_4=(x_0,0.75)^T-(\alpha_0,-\beta_0)^T,
\end{aligned}
\end{equation}
with $0<\beta_0,\alpha_0<1$ such that $\bvec{Q}_0\ne \bvec{0}$, then the analytical solution of \eqref{newODE} with \eqref{initODE} is given as
\begin{equation}\label{eq:symmetry of sol4c}
\begin{aligned}
&\va_1(t)=(x_0,0.25)^T+(\alpha(t),\beta(t))^T-2t\bvec{Q}_0, \quad \va_2(t)=(x_0,0.25)^T-(\alpha(t),\beta(t))^T-2t\bvec{Q}_0,\\ &\va_3(t)=(x_0,0.75)^T+(\alpha(t),-\beta(t))^T-2t\bvec{Q}_0, \quad \va_4(t)=(x_0,0.75)^T-(\alpha(t),-\beta(t))^T-2t\bvec{Q}_0,
\end{aligned}
\end{equation}
where ($\alpha(t)$,$\beta(t)$) is the solution of 
\begin{equation}\label{eq:dynamics of alphabeta4c}
\left\{
\begin{aligned}
  &\dot{\alpha}=2\left(\p_yF(2\alpha,2\beta)-\p_y(0,2\beta-0.5) \right),\\
  &\dot{\beta}=2(-\p_xF(2\alpha,2\beta)+\p_xF(2\alpha,0.5)),
  \end{aligned}
  \right.
\end{equation}
with the initial data \eqref{eq:reduced initial}.
\end{Lem} }

\begin{figure}[htp!]
\begin{center}
\begin{tabular}{cccc}
\includegraphics[height=4.5cm]{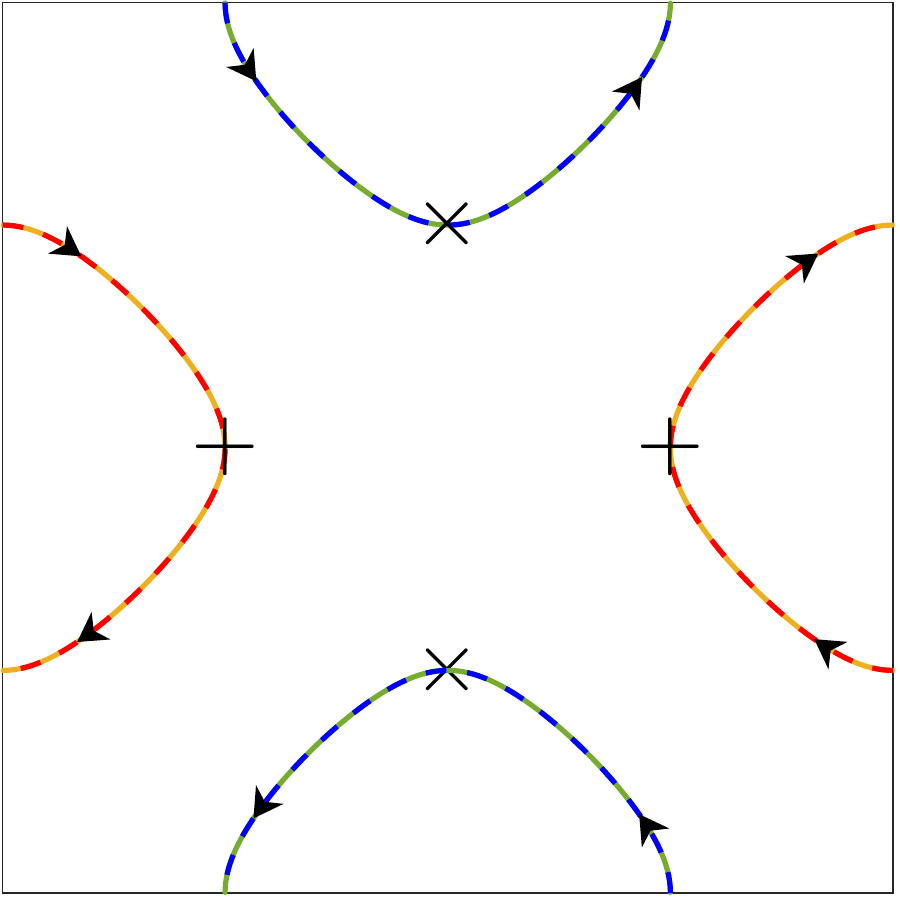}
&\includegraphics[width=4.5cm]{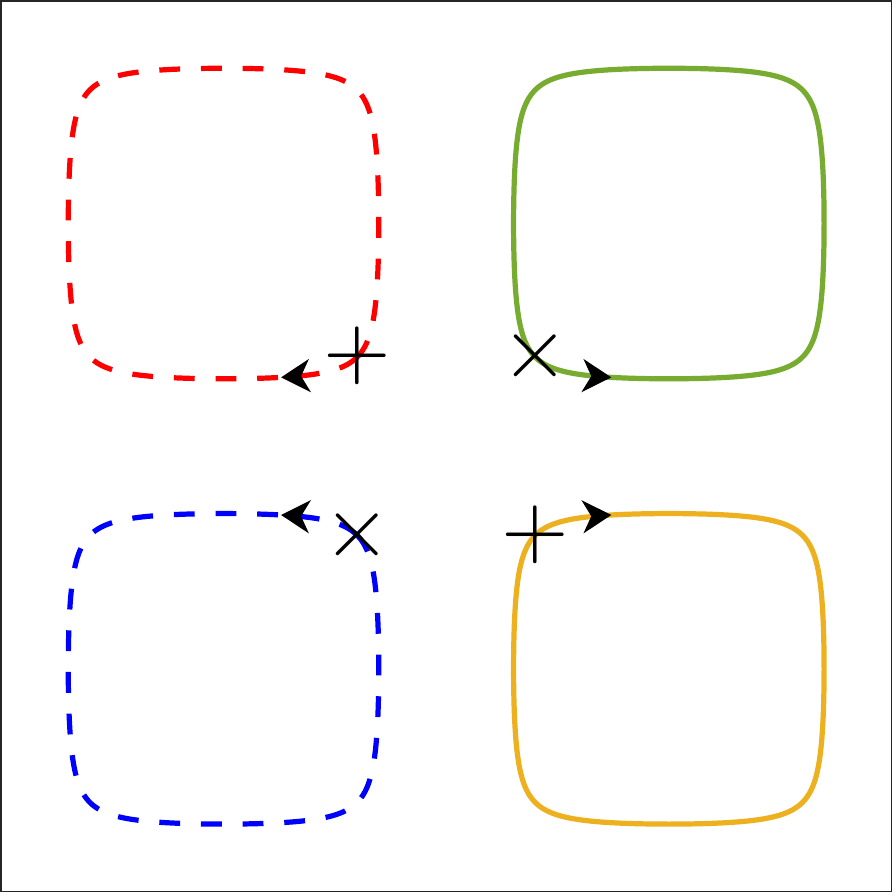}
&\includegraphics[height=4.5cm]{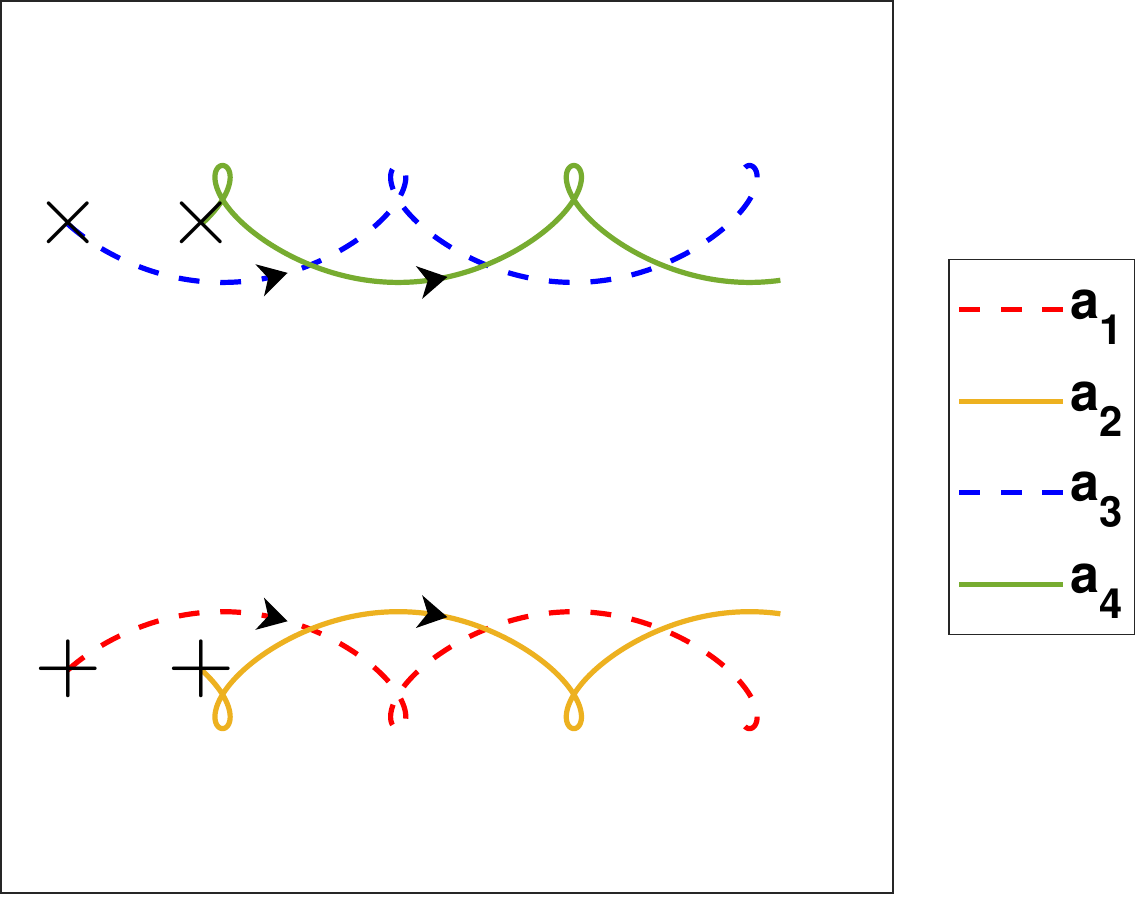}
\end{tabular}
\end{center}
\caption{Trajectories of the reduced dynamics law \eqref{newODE} with $N=2$ for 
different initial datum: (i) in \eqref{eq:initial of so} with $\alpha_0=-0.25$ and
$\beta_0=0$ (left), (ii) in \eqref{initN21} with $\alpha_0=-0.1$ and $\beta_0=0.1$
(middle), and (iii) in \eqref{eq:initial of so4c} with $x_0=0.15$, $\alpha_0=-0.075$ and $\beta_0=0$ (right).}\label{fig:aaa}
\end{figure}

\Bluemark{To illustrate the solution of \eqref{newODE} with the initial data \eqref{initODE} for a few special setups, 
we solve it numerically by adopting the fourth-order Runge-Kutta method with time step $\Delta t=10^{-4}$. 
Figure \ref{fig:2vpaths} plots the solutions of \eqref{newODE} with $N=1$ and different initial datum in 
\eqref{initODE}, i.e. different $\bvec{Q}_0$, to illustrate
the dynamics described in Lemma \eqref{ASODEN1}; and Figure \ref{fig:aaa} shows 
the solutions of \eqref{newODE} with $N=2$ and different initial datum in 
\eqref{eq:initial of so} to illustrate
the dynamics described in Lemma \ref{lem:2circles}, in  \eqref{initN21} to illustrate
the dynamics described in Lemma \ref{SODEN21}, and in \eqref{eq:initial of so4c} to illustrate
the dynamics described in Lemma \ref{SODEN22}.}

\section{Conclusion}\label{sec:conclusion}

\blackmark{A new reduced dynamical law} for quantized vortex dynamics of the nonlinear Sch\"{o}dinger equation (NLSE) on \blackmark{the torus} with non-vanishing momentum was established
when the vortex core size $\varepsilon\to 0$.
It is governed by a Hamiltonian flow driven by a renormalized energy on \blackmark{the torus}
and it collapses to \blackmark{the reduced dynamical law} obtained in \cite{RN8} for NLSE
on \blackmark{the torus} with vanishing momentum. The key step is to adopt  a new canonical harmonic map on \blackmark{the torus} to include the effect of the non-vanishing momentum  into the dynamics. Extension of \blackmark{the reduced dynamical law} for NLSE
on \blackmark{torus} with arbitrary length and width was discussed. Finally,
three first integrals of \blackmark{the reduced dynamical law} were presented and analytical
solutions were obtained with several initial setups with symmetry.

\section*{CRediT authorship contribution statement}
\textbf{Weizhu Bao:} Conceptualization, Validation, Supervision, Writing - review \& editing.  \textbf{Huaiyu Jian:} Conceptualization, Validation, Supervision, Writing - review \& editing. \textbf{Yongxing Zhu:} Conceptualization, Methodology, Writing - original draft, Writing - review \& editing.

\section*{Data availability}
No data was used for the research described in the article.

\section*{Acknowledgments}
This work was partially supported by the China Scholarship Council (Y. Zhu),
by the Ministry of Education of Singapore under its  AcRF Tier 2 funding 
MOE-T2EP20122-0002 (A-8000962-00-00) (W. Bao), and by
the National Natural Science Foundation of China grant 12141103
(H. Jian). Part of the work was done when the first two authors were visiting
the Institute for Mathematical Science in 2023.
The referees are thanked for their valuable
comments, which improved the quality of the paper.

\section*{Declaration of competing interest}
The authors declare that they have no known competing financial interests or personal relationships that could have appeared to influence the work reported in this paper.

\bibliographystyle{model1-num-names}

\end{document}